\documentclass[12pt,a4paper]{article}
\usepackage[T1]{fontenc}
\usepackage{multirow,scalerel}
\usepackage{setspace}
\usepackage[linesnumbered,ruled,vlined]{algorithm2e}

\usepackage{algpseudocode}
\usepackage[numbered]{bookmark}

\usepackage{appendix}

\usepackage{amsfonts,enumerate,amsmath}
\usepackage{latexsym}
\usepackage{color}
\usepackage{natbib}
\usepackage{array}

\usepackage{verbatim}
\usepackage{graphicx}
\usepackage{array}
\usepackage{float,bm}
\usepackage{multirow}

\renewcommand{\theequation}{\arabic{section}.\arabic{equation}}

\newcommand{\qed}{\hfill\Box}

\newtheorem{assumps}{Assumption}

\graphicspath{{./images/}}

\renewcommand{\cite}{\citet*}

\newcommand{\bqn}{\begin{eqnarray*}}
\newcommand{\eqn}{\end{eqnarray*}}

\DeclareMathOperator{\tr}{tr}         
      \DeclareMathOperator{\diag}{diag}

\numberwithin{equation}{section}  

\newtheorem{thm}{Theorem}
\newtheorem{prop}{Proposition}
\newtheorem{lem}{Lemma}

\hypersetup{
  colorlinks   = true, 
  urlcolor     = blue, 
  linkcolor    = blue, 
  citecolor   = blue 
}

\usepackage{array}
\usepackage{multirow}
\usepackage{bibunits}
\defaultbibliography{Cov_est_elliptical-aos_template}
\defaultbibliographystyle{apalike}

\begin{document}


\title{
Sub-Gaussian High-Dimensional Covariance Matrix Estimation under Elliptical Factor Model with $2+\varepsilon$th Moment
}

\author{Yi Ding \thanks{Faculty of Business Administration, University of Macau, Macau. Research is supported in part by NSFC 72101226 of China, and SRG2023-00001-FBA, MYRG-GRG2023-00119-FBA  of the University of Macau.  Email: yiding@um.edu.mo}
 \and Xinghua Zheng \thanks{Department of ISOM, Hong Kong University of Science and Technology, Clear Water Bay, Kowloon, Hong Kong. Research is supported in part by RGC grants GRF-16304521 and GRF~16304019 of the HKSAR. Email: xhzheng@ust.hk}}
\date{\today}

\maketitle
\begin{abstract}
\noindent 
We study the estimation of high-dimensional covariance matrices under elliptical factor models with $2+\varepsilon$th moment. 
For such heavy-tailed data, robust estimators like the Huber-type estimator in \cite{fan2018large} can not achieve sub-Gaussian convergence rate. 
In this paper, we develop an idiosyncratic-projected self-normalization (IPSN) method to remove the effect of  heavy-tailed scalar parameter, 
and propose a robust pilot estimator for the scatter matrix that achieves the  sub-Gaussian rate. We further develop an estimator of the 
covariance matrix and show that it achieves a faster convergence rate than the generic POET estimator in  \cite{fan2018large}.

 \vskip 0.3cm



\end{abstract}

\section{Introduction}

\subsection{Covariance Matrix Estimation for Heavy-tailed Distributions}

There is an extensive interest in  robust estimation of the covariance matrices for heavy-tailed data, which are commonly seen in various fields such as finance, economics and biology. Huber-type robust estimators (\cite{minsker2018sub,minsker2020robust}), and truncation or shrinkage methods (\cite{fan2021shrinkage}) are popular approaches in this area. Notably, these studies typically assume  that data have finite~4th moment.

\cite{avella2018robust}  study the high-dimensional covariance matrix estimation for cases where data have finite 4th or~\mbox{$2+\varepsilon$th} moment with~\mbox{$\varepsilon<2$}. Under sparsity assumptions on the covariance matrix or its inverse matrix, their approach involves applying thresholding techniques (\cite{cai2011adaptive,cai2016estimating}) on a robust pilot matrix. The pilot matrix combines a Huber-type M-estimator of variances and a rank-based robust correlation matrix estimator such as the marginal Kendall's tau (\cite{liu2012high,xue2012regularized}). Under the high-dimensional setting where the dimension  $p$ and the sample size  $n$ both go to infinity,  if the data have 4th moment, \cite{avella2018robust} show that their estimator achieves the convergence rate of $\sqrt{(\log p)/n}$, which is the same convergence rate as in the sub-Gaussian case. However, for the case when only  $2+\varepsilon$th moment exists, the convergence rate reduces to $((\log p)/n)^{\varepsilon/(2+\varepsilon)}$; see Proposition 6 therein.

\subsection{Elliptical Factor Models}
Due to their flexibility in modeling heavy-tailedness and capturing tail-dependences, elliptical models are widely used to model financial returns (\cite{ZL11,han2014scale}), gene microarrays (\cite{avella2018robust}) and brain images (\cite{han2018eca}).
An elliptical model is typically defined as follows:
\begin{equation}\label{ellip_model}
\mathbf{y}_t=\boldsymbol{\mu}+\sqrt{p}\xi_t\boldsymbol{\Sigma}_0^{1/2}\mathbf{u}_t,
\end{equation}
where $\boldsymbol{\mu}=(\mu_1, ..., \mu_p)^{T}$ is the expectation, $\xi_t\geq 0$ is a  random scalar variable, $\boldsymbol{\Sigma}_0$ is a $p\times p$ nonnegative matrix, and $\mathbf{u}_t$ is a length-$p$ vector uniformly distributed on the unit sphere $\mathcal{S}^{p-1}$ and independent of $\xi_t$.
The covariance matrix of $(\mathbf{y}_t)$ is $\boldsymbol{\Sigma}=E(\xi_t^2)\boldsymbol{\Sigma}_0$. Usually, $\boldsymbol{\Sigma}_0$ is  referred to as the scatter matrix (\cite{fan2018large}). For identification purposes, $\boldsymbol{\Sigma}_0$ can be normalized to have a trace of $p$. Popular estimators of the scatter matrix include the minimum covariance determinant (MCD) and minimum volume ellipsoid (MVE) estimators (\cite{rousseeuw1984least,rousseeuw1985multivariate}), Maronna's M-estimator (\cite{maronna1976robust}) and
Tyler's M-estimator (\cite{tyler1987distribution}). However, these estimators can not handle the case when $p> n$; see, for example, the discussions in \cite{van2009minimum}, \cite{zhang2016robust,zhang2016marvcenko} and \cite{hubert2018minimum}.

High-dimensional data such as financial returns usually exhibit a factor structure. 
\cite{FFL08, FLM11, fan2013large} utilize  factor models coupled with idiosyncratic sparsity to estimate the high-dimensional covariance matrix under such settings. These studies focus on the case where  data are sub-Gaussian or sub-Exponential.

\cite{fan2018large} consider elliptical factor models under which only $4$th moment exists and propose a generic Principal Orthogonal ComplEment Thresholding (POET)  procedure to estimate the high-dimensional covariance matrix. Their approach relies on three components,~$\widehat{\boldsymbol{\Sigma}}$,~$\widehat{\boldsymbol{\Lambda}}_K$  and~$\widehat{\boldsymbol{\Gamma}}_K$, which  are robust pilot estimators of the covariance matrix $\boldsymbol{\Sigma}$,  its leading eigenvalues $\boldsymbol{\Lambda}_{K}$ and eigenvectors $\boldsymbol{\Gamma}_{K}$, respectively. The theoretical properties of the generic POET procedure crucially depend on the consistency of the pilot components, which have sub-Gaussian convergence rate under the $4$th moment condition:
\begin{align}
&\|\widehat{\boldsymbol{\Sigma}}-\boldsymbol{\Sigma}\|_{\max}=O_p\Bigg(\sqrt{\frac{\log p}{n}}\Bigg),\label{fanrate1}\\
&\|\widehat{\boldsymbol{\Lambda}}_K{\boldsymbol{\Lambda}}^{-1}_K-\mathbf{I}\|_{\max}=O_p\Bigg(\sqrt{\frac{\log p}{n}}\Bigg),\quad \text{and}\label{fanrate2}\\
&\|\widehat{\boldsymbol{\Gamma}}_K-{\boldsymbol{\Gamma}}_K\|_{\max}=O_p\Bigg(\sqrt{\frac{\log p}{np}}\Bigg).\label{fanrate3}
\end{align}
When only  $2+\varepsilon$th moment exists, the convergence rates in equations \eqref{fanrate1} and \eqref{fanrate2} reduce to $((\log p)/n)^{\varepsilon/(2+\varepsilon)}$.

\subsection{Our Contributions}

Infinite kurtosis and cross-sectional tail dependence are commonly seen in financial data. Stock returns contain large and dependent jumps (e.g., \cite{ding2023stock}), which  have slowly decaying  tails that suggest infinite 4th moment (\cite{bollerslev2013jump}). For  data with only $2+\varepsilon$th moment,  the aforementioned covariance matrix estimation methods  can not achieve the sub-Gaussian rate.
An important and challenging question is:
\\\\
\emph{Is there a covariance matrix estimator that can achieve a higher rate than $((\log p)/n)^{\varepsilon/(2+\varepsilon)}$ or even  the sub-Gaussian rate of $\sqrt{(\log p)/n}$ when only $2+\varepsilon$th moment exists?}
\\\\
We address this question in this paper. We study the robust estimation of high-dimensional covariance matrices under elliptical factor models with only $2+\varepsilon$th moment.
It is worth mentioning that in the elliptical model \eqref{ellip_model},  the scatter matrix~$\boldsymbol{\Sigma}_0$ determines the cross-sectional dependence  and plays a central role in various applications. For example, the eigenstructure of $\boldsymbol{\Sigma}_0$ governs the principal component analysis (PCA) of $(\mathbf{y}_t)$. In financial applications, the inverse of the scatter matrix,~$\boldsymbol{\Sigma}_0^{-1}$, is used in constructing the minimum variance portfolio.  The scalar component, $E(\xi_t^2)$, on the other hand, is not relevant in the above applications.
These observations motivate us to estimate the two components, $E(\xi^2_t)$ and $\boldsymbol{\Sigma}_0$, separately.

In this paper, we propose an estimator of the scatter matrix $\boldsymbol{\Sigma}_0$ that achieves the  sub-Gaussian convergence rate under the condition that $(\mathbf{y}_t)$  has only $2+\varepsilon$th moment. We summarize the construction of
our estimator below.

First,
we estimate $(\xi_t)$ and remove the effect of $(\xi_t)$ in $(\mathbf{y}_t)$. \cite{ZL11} develop a self-normalization approach to remove $(\xi_t)$ in elliptical models. Specifically, when~$p$ is large,~$(\xi_t)$ can be consistently estimated  by $(\|\mathring{\mathbf{y}}_t\|_2/\sqrt{p})$, where $(\mathring{\mathbf{y}}_t)$ is a properly demeaned~$(\mathbf{y}_t)$. The consistency of $(\|\mathring{\mathbf{y}}_t\|_2/\sqrt{p})$  in estimating $(\xi_t)$ requires that the data do not admit a strong factor structure. Under the elliptical factor model, such a consistency no longer holds. In this paper, we develop an idiosyncratic projected self-normalization (IPSN) approach that allows the existence of factors.

Next, we use the idiosyncratic-projected self-normalized data to construct the pilot components for the generic POET estimation of $\boldsymbol{\Sigma}_0$. We  use the sample covariance matrix of the idiosyncratic-projected self-normalized data and standardize it to have trace $p$. The resulting estimator is the pilot covariance matrix estimator, $\widehat{\boldsymbol{\Sigma}}_0$. The pilot leading eigenvalues and leading eigenvector estimators, $\widehat{\boldsymbol{\Lambda}}_{0K}$ and $\widehat{\boldsymbol{\Gamma}}_{0K}$, are simply the leading eigenvalues and eigenvectors of $\widehat{\boldsymbol{\Sigma}}_0$. We show that our pilot estimators achieve the sub-Gaussian rate under only~\mbox{2+$\varepsilon$}th moment condition. We then apply the generic POET procedure on the proposed~IPSN pilot estimators to  obtain the scatter matrix estimator, denoted by GPOET-IPSN. We show that the  GPOET-IPSN estimator achieves the sub-Gaussian convergence rate under a conditional sparsity assumption.

Finally, to estimate the covariance matrix $\boldsymbol{\Sigma}$, we  construct a consistent estimator of $E(\xi_t^2)$ based on  $(\|\mathring{\mathbf{y}}_t\|_2^2)$. We show that the  estimator  achieves the same converge rate as the oracle estimator based on the unobserved $(\xi_t)$ process. Combining our proposed \mbox{GPOET-IPSN estimator} of~$\boldsymbol{\Sigma}_0$ with the estimator of $E(\xi_t^2)$ yields our estimator of $\boldsymbol{\Sigma}$, which we show  achieves a faster convergence rate than the  generic POET estimator in \cite{fan2018large}.


The favorable properties of our proposed  estimator are clearly demonstrated in the numerical studies. Simulation studies show that the \mbox{GPOET-IPSN} estimator performs robustly well under various heavy-tailedness settings and achieves a higher estimation accuracy than the robust generic POET estimator in \cite{fan2018large} and the POET estimator in \cite{fan2013large}.
Empirically, we apply our proposed estimator to construct  minimum variance portfolios using S\&P~500 Index constituent stocks. We find that our portfolio has a statistically significantly lower out-of-sample risk than benchmark portfolios.





 The rest of this paper is organized as follows. We present our theoretical results in Section~\ref{main}. Simulations and empirical studies are presented in Sections~\ref{Simu} and~\ref{Empi}, respectively.  We conclude in Section~\ref{Conc}. The proofs of the theorems and propositions are presented in Section~\ref{Proof}. The proofs of the lemmas are collected in the Supplementary Material (\cite{DZ24_supp}).

 We use the following notation throughout the paper.
For any matrix~$\mathbf{A}=(A_{ij})$, its spectral norm is defined as~$\|\mathbf{A}\|_2=\max_{\|\mathbf{x}\|_2\leq 1}\sqrt{\mathbf{x}^T\mathbf{A}^T\mathbf{A}\mathbf{x}}$, where $\|\mathbf{x}\|_2=\sqrt{\sum x_i^2}$ for any vector~$\mathbf{x}=(x_i)$; the Frobenius norm is defined as~$\|\mathbf{A}\|_{F}=\sqrt{\sum_{i,j}A_{ij}^2}$; and the relative Frobenius norm is defined as $\|\mathbf{A}\|_{\boldsymbol{\Sigma}}=\|\boldsymbol{\Sigma}^{-1/2}\mathbf{A}\boldsymbol{\Sigma}^{-1/2}\|_F/\sqrt{p}$, where $\boldsymbol{\Sigma}$ is a~\mbox{$p\times p$} positive-definite matrix.

\section{Main Results}\label{main}
\subsection{Settings and Assumptions}
Suppose that $(\mathbf{y}_t)$ follows the elliptical model \eqref{ellip_model}. We make  the following assumptions.
\begin{assumps}\label{asump1} $(\xi_t)$ are \mbox{i.i.d.} and independent with $(\mathbf{u}_t)$. Moreover, there exist constants $0<c<C<\infty$ such that $E(\xi_t^{2+\varepsilon})<C$ for some $0<\varepsilon< 2$ and $P(|\xi_t|>c)=1$.
\end{assumps}
\begin{assumps}\label{sigma_bound} There exist constants $0<c<C<\infty$ such that  $c<\min_{1\leq 1 \leq p}(\boldsymbol{\Sigma}_0)_{ii}\leq\max_{1\leq 1 \leq p}(\boldsymbol{\Sigma}_0)_{ii}<C$.
\end{assumps}

Write the eigen decomposition of $\boldsymbol{\Sigma}_0$ as $\boldsymbol{\Sigma}_0=\mathbf{\Gamma}\boldsymbol{\Lambda}_0\mathbf{\Gamma}^T$, where  $\boldsymbol{\Lambda}_0=\diag(\lambda_{0;1}, \lambda_{0;2}, ..., \lambda_{0;p})$, $\lambda_{0;1}\geq \lambda_{0;2}...\geq \lambda_{0;p}$  are the eigenvalues of $\boldsymbol{\Sigma}_0$, and $\mathbf{\Gamma}=(\phi_1,...,\phi_p)$ is the corresponding matrix of eigenvectors. We assume that factors exist so that the following eigen-gap condition holds.
\begin{assumps}\label{assump_factor_structure}There exist a $K\in \mathbb{N}$ and $\delta>0$ such that  $\delta p\leq \lambda_{0;i}\leq p/\delta$, and $0<\delta \leq\lambda_{0;i}\leq 1/\delta$ for $i=K+1, ..., p$.
\end{assumps}
{Widely used approaches for estimating the number of factors such as those developed in \cite{bai2002determining}, \cite{onatski2010determining} and \cite{ahn2013eigenvalue} require the existence of~4th moment. Under the elliptical factor model with $2+\varepsilon$th moment, the number of factors~$K$ can be consistently estimated using the robust estimator in \cite{yu2019robust}, which utilizes eigenvalue ratios of spatial Kendall's tau matrices. In the following analysis, we assume that the number of factors is known.
}

\subsection{Idiosyncratic-Projected Self-Normalized Estimator of Covariance Matrix}\label{approach}

Observe that under the model \eqref{ellip_model}, we can write $\mathbf{y}_t$ as follows:
\begin{equation}\label{Y:model}
\mathbf{y}_t=\boldsymbol{\mu}+\frac{\sqrt{p}\xi_t}{\|\mathbf{z}_t\|_2}\boldsymbol{\Sigma}^{1/2}_0\mathbf{z}_t,
\end{equation}
where $\mathbf{z}_t\sim N(0, \mathbf{I}_p)$, and $\xi_t$ and $\mathbf{z}_t$ are independent. The heavy-tailedness of $(\mathbf{y}_t)$ is solely determined by $(\xi_t)$. Motivated by this observation, we separately estimate  $E(\xi^2_t)$ and $\boldsymbol{\Sigma}_0$. To do so, we develop an idiosyncratic-projected self-normalization (IPSN) approach.
We explain the construction of our proposed estimator below.\\\\
{\bf Step I. Estimate $\boldsymbol{\mu}$}

Given that $(\mathbf{y}_t)$ is heavy-tailed, in Step I, we estimate $\boldsymbol{\mu}$ using a Huber estimator. Specifically, the estimator $\widehat{\boldsymbol{\mu}}=(\widehat{\mu}_i)_{1\leq i\leq p}$ solves
\begin{equation}\label{hat:mu}
\sum_{t=1}^n\Psi_H(y_{i,t}-\widehat{\mu}_i)=0\quad\text{ for} \quad i=1, ..., p,
\end{equation}
where $\Psi_H$ is the Huber function:
\begin{equation}\label{Psi_H}
\Psi_H(x)=\min\Big(H, \max (-H, x)\Big)
\end{equation}
with $H$  a tuning parameter, which goes to infinity with $n$ and satisfies $H=O(\sqrt{n/(\log p)})$. Numerically, the tuning parameter $H$ can be chosen with  cross-validation (\cite{sun2020adaptive}). \\\\

\noindent{\bf Step II. Construct idiosyncratic-projected self-normalized {variables} $(\widehat{\mathbf{X}}_t)$}

In Step II, we normalize observations  to remove the effect of $(\xi_t)$. We use  the $\ell_2$ norm of  estimated idiosyncratic variables to conduct normalization. Specifically,  we first compute the spatial Kendall's tau matrix:
$$
\widehat{\boldsymbol{\Sigma}}_2=\frac{2}{n(n-1)}\sum_{i<i'}\frac{(\mathbf{y}_i-\mathbf{y}_{i'})(\mathbf{y}_i-\mathbf{y}_{i'})^T}{\|\mathbf{y}_i-\mathbf{y}_{i'}\|_2^2},
$$
 get its leading eigenvectors, denoted by $\widehat{\boldsymbol{\Gamma}_K}_{ED}$, and   use them as a robust estimator of  $\boldsymbol{\Gamma}_{K}=(\phi_1,...,\phi_K)$.  Then, we
define a $(p-K)\times p$ projection matrix $\widehat{\mathbf{P}}_I$ as follows:
$$\widehat{\mathbf{P}}_I^T=\text{\rm Null}(\widehat{\boldsymbol{\Gamma}_K}_{ED}),$$
that is, $\widehat{\mathbf{P}}_I^T$ is the $p\times (p-K)$  matrix that spans the null space of $\widehat{\boldsymbol{\Gamma}_K}_{ED}$ and satisfies \mbox{$\widehat{\mathbf{P}}_I\widehat{\boldsymbol{\Gamma}_K}_{ED}=\mathbf{0}$} and $\widehat{\mathbf{P}}_I\widehat{\mathbf{P}}_I^T=\mathbf{I}_{p-K}$.

The goal of the idiosyncratic projection is to remove the strong cross-sectional dependence in the data, thereby enabling consistent estimation of the {scalar component} $(\xi_t)$. This step is crucial. As can be seen in Section \ref{simu_result} below,  when there is strong cross-sectional dependence, the self-normalization approach in \cite{ZL11} no longer works in estimating the time-varying scalar component~$(\xi_t)$.

Next, we define the idiosyncratic projected self-normalized variables $(\widehat{\mathbf{X}}_t)$:
\begin{equation}\label{x:sn:case2}
\widehat{\mathbf{X}}_t=\sqrt{p}\frac{\mathbf{y}_{t}-\widehat{\boldsymbol{\mu}}}{\|\widehat{\mathbf{P}}_I(\mathbf{y}_{t}-\widehat{\boldsymbol{\mu}})\|_2}, \quad\quad t=1, ..., n.
\end{equation}
The normalization terms $(\|\widehat{\mathbf{P}}_I(\mathbf{y}_{t}-\widehat{\boldsymbol{\mu}})\|_2)$ are approximately proportional to $(\xi_t)$, hence~$(\widehat{\mathbf{X}}_t)$ are approximately proportional to $(\boldsymbol{\Sigma}_0^{1/2}\mathbf{u}_t)$. \\\\

\noindent{\bf Step III. Construct pilot estimators of $\boldsymbol{\Sigma}_0$}

In Step III, we estimate $\boldsymbol{\Sigma}_0$ with the following estimator:
\begin{equation}\label{hatSig0}
\widehat{\boldsymbol{\Sigma}}_0=\frac{\widehat{\eta}}{n}\sum_{t=1}^n \widehat{\mathbf{X}}_t\widehat{\mathbf{X}}_t^T,
\end{equation}
where
\begin{equation}\label{gam_norm}
\widehat{\eta}= \frac{p}{\tr(\sum_{t=1}^n \widehat{\mathbf{X}}_t\widehat{\mathbf{X}}_t^T/n)}.
\end{equation}
By doing so, we essentially normalize $\mathbf{y}_{t}-\widehat{\boldsymbol{\mu}}$ with
\begin{equation}\label{hatxi_t}
\widehat{\xi}_t=\frac{\|\widehat{\mathbf{P}}_I(\mathbf{y}_{t}-\widehat{\boldsymbol{\mu}})\|_2}{\sqrt{\widehat{\eta}}}.
\end{equation}
By the definition of $\widehat{\eta}$ in equation \eqref{gam_norm},  we ensure that $\tr(\widehat{\boldsymbol{\Sigma}}_0)=p$.

Write the eigen-decomposition of $\widehat{\boldsymbol{\Sigma}}_0$ as $\widehat{\boldsymbol{\Sigma}}_0=\widehat{\mathbf{\Gamma}}\widehat{\boldsymbol{\Lambda}}_0\widehat{\mathbf{\Gamma}}^T$, where $\widehat{\boldsymbol{\Lambda}}_0=\diag(\widehat{\lambda}_{0;1},...,\widehat{\lambda}_{0;p})$ is the matrix of eigenvalues with $\widehat{\lambda}_{0;1}\geq \widehat{\lambda}_{0;2}\geq...\geq \widehat{\lambda}_{0;p}$, and $\widehat{\mathbf{\Gamma}}=(\widehat{\phi}_1,...,\widehat{\phi}_p)$ is the matrix of eigenvectors. We then use the leading eigenvectors  and leading eigenvalues of~$\widehat{\boldsymbol{\Sigma}}_0$ as the corresponding pilot estimators: $$\widehat{\mathbf{\Gamma}}_{K}=(\widehat{\phi}_{1},...,\widehat{\phi}_K), \quad\text{ and }\quad \widehat{\boldsymbol{\Lambda}}_{0K}=\diag(\widehat{\lambda}_{0;1},...,\widehat{\lambda}_{0;K}).$$
The  estimators $\widehat{\boldsymbol{\Sigma}}_0$, $\widehat{\boldsymbol{\Gamma}}_{K}$ and $\widehat{\boldsymbol{\Lambda}}_{0K}$ will be used in the generic POET procedure for estimating~$\boldsymbol{\Sigma}_0$.

In the generic POET procedure  proposed in \cite{fan2018large},  the leading eigenvectors are estimated using the spatial Kendall's tau, and pilot components of the covariance matrix and its leading eigenvalues are from a Huber estimator of the covariance matrix.
As these estimators are from different sources, the resulting pilot idiosyncratic covariance matrix  estimator is not guaranteed to be semi-positive definite. Our approach, on the other hand, does not have such an issue because our pilot estimators of the leading eigenvalues and leading eigenvectors are from the eigen decomposition of a same covariance matrix estimator.
\\\\
\noindent{\bf Step IV. Estimate $E(\xi_t^2)$}

In order to estimate~$\boldsymbol{\Sigma}$,  we also need to estimate $E(\xi_t^2)$.
We propose the following robust estimator $\widehat{E(\xi_t^2)}$ that solves
\begin{equation}\label{xi2_est}
\sum_{i=1}^n \Psi_H\Bigg(\frac{\|\mathbf{y}_t-\widehat{\boldsymbol{\mu}}\|_2^2}{p}-\widehat{E(\xi_t^2)}\Bigg)=0,
\end{equation}
where $H\asymp n^{\min(1/(1+\varepsilon/2),1/2)}$. Again, the tuning parameter $H$ can be chosen  using cross-validation.
\\\\
\noindent{\bf Step V. Generic POET estimators of $\boldsymbol{\Sigma}_0$ and $\boldsymbol{\Sigma}$}

Finally, we  estimate the scatter matrix $\boldsymbol{\Sigma}_0$ and the covariance matrix $\boldsymbol{\Sigma}$.
To estimate~$\boldsymbol{\Sigma}_0$, we apply the generic POET procedure  on our proposed  idiosyncratic-projected self-normalized pilot estimators, $\widehat{\boldsymbol{\Sigma}}_0$,  $\widehat{\boldsymbol{\Lambda}}_{0K}$ and $\widehat{\boldsymbol{\Gamma}}_{K}$  from Step III. Specifically,
 the estimator of $\boldsymbol{\Sigma}_0$ is \begin{equation}\label{gPOET}
\widehat{\boldsymbol{\Sigma}}^{\tau}_0=\widehat{\boldsymbol{\Gamma}}_{K}\widehat{\boldsymbol{\Lambda}}_{0K}\widehat{\boldsymbol{\Gamma}}_{K}^T+\widehat{\boldsymbol{\Sigma}}_{0u}^{\tau},
\end{equation}
where $\widehat{\boldsymbol{\Sigma}}_{0u}^{\tau}$ is obtained by applying the adaptive thresholding method (\cite{cai2011adaptive}) to $\widehat{\boldsymbol{\Sigma}}_{0}-\widehat{\boldsymbol{\Gamma}}_{K}\widehat{\boldsymbol{\Lambda}}_{0K}\widehat{\boldsymbol{\Gamma}}_{K}^T$.

To estimate $\boldsymbol{\Sigma}$, we combine $\widehat{E(\xi_t^2)}$ from Step IV with $\widehat{\boldsymbol{\Sigma}}^{\tau}_0$ in equation \eqref{gPOET}:
\begin{equation}\label{gpoet_s}
\widehat{\boldsymbol{\Sigma}}^{\tau}=\widehat{E(\xi_t^2)}\, \widehat{\boldsymbol{\Sigma}}^{\tau}_0.
\end{equation}
The corresponding idiosyncratic covariance matrix estimator is $\widehat{\boldsymbol{\Sigma}}_{u}^{\tau}=\widehat{E(\xi_t^2)}\,\widehat{\boldsymbol{\Sigma}}_{0u}^{\tau}$.
It can be seen that $\widehat{\boldsymbol{\Sigma}}^{\tau}$ is equivalent to  the generic POET estimator using the following pilot estimators: $\widehat{\boldsymbol{\Sigma}}=:\widehat{E(\xi_t^2)}\,\widehat{\boldsymbol{\Sigma}}_0$,  $\widehat{\boldsymbol{\Lambda}}_{K}=:\widehat{E(\xi_t^2)}\,\widehat{\boldsymbol{\Lambda}}_{0K}$, and $\widehat{\boldsymbol{\Gamma}}_{K}$. \\\\
We summarize our IPSN approach in the following algorithm. \\\\
\begin{tabular}{p{0.001\textwidth}p{0.1\textwidth}p{0.8\textwidth}}
\hline
\multicolumn{2}{l}{{\bf Algorithm}:}& Idiosyncratic-Projected Self-Normalization (IPSN) for Covariance  Matrix Estimation \\\hline
&{\bf Input:}& {$(\mathbf{y}_t)_{t\leq n}$, $K$.}\\
&{\bf Output:}& {$\widehat{\boldsymbol{\Sigma}}^{\tau}_0$, $\widehat{\boldsymbol{\Sigma}}^{\tau}$.}\\
& Step I.& Compute $\widehat{\boldsymbol{\mu}}$ via \eqref{hat:mu}.\\
&Step II.& Compute $(\widehat{\mathbf{X}}_t)_{t\leq n}$ via \eqref{hat:mu}.\\
&Step III.& Compute $\widehat{\eta}$ and $(\widehat{\xi}_t)_{t\leq n}$ via \eqref{gam_norm} and \eqref{hatxi_t}, respectively. Compute  $\widehat{\boldsymbol{\Sigma}}_0$ via \eqref{hatSig0},  its first $K$ eigenvalues $\widehat{\boldsymbol{\Lambda}}_{0K}$, and the corresponding eigenvectors $\widehat{\mathbf{\Gamma}}_{K}$. \\
&Step IV.& Compute $\widehat{E(\xi_t^2)}$  via \eqref{xi2_est}.\\
&Step V.& Compute $\widehat{\boldsymbol{\Sigma}}^\tau_0$ via \eqref{gPOET}, and $\widehat{\boldsymbol{\Sigma}}^\tau$ via \eqref{gpoet_s}.\\
\hline
\end{tabular}

\subsection{Theoretical Properties}

In this subsection, we present the asymptotic properties of our proposed estimators.

\subsubsection{Convergence Rates of Idiosyncratic-Projected Self-Normalized Pilot Estimators}

We start with the properties of the idiosyncratic-projected self-normalization. The key to the success of our approach lies in  approximating  $(\xi_t)$ with~$(\widehat{\xi}_t)$ defined in equation~\eqref{hatxi_t}. Intuitively, if $(\widehat{\xi}_t)$ is close to $(\xi_t)$, then $((\mathbf{y}_{t}-\widehat{\boldsymbol{\mu}})/\widehat{\xi}_t)$ is close to $\sqrt{p}(\boldsymbol{\Sigma}_0^{1/2}\mathbf{u}_t)$, which is absent of the heavy-tailed scalar process $(\xi_t)$.
The next proposition gives the properties of $~(\widehat{\xi}_t)$.
\begin{prop}\label{lem_denominator_xi}Under Assumptions \ref{asump1}--\ref{assump_factor_structure}, if in addition,  $p, n\to \infty$ and satisfy \\\mbox{$(\log(p))^{2+\gamma}=o(n)$} for some $\gamma>0$ and $\log(n)=O(\log (p))$, then there exists a constant $c_1>0$ such that
\begin{equation}\label{prob_bound_1_xi}
P\Bigg(\min_{1\leq t\leq n}\widehat{\xi}_t>c_1\Bigg)\to 1,
\end{equation}
and
\begin{equation}\label{sum_bound_1_xi}
\frac{1}{n}\sum_{t=1}^n\Big(\frac{\xi_t^2}{\widehat{\xi}_t^2}-1\Big)^2=O_p\Big(\frac{\log p}{p}+\frac{\log p}{n}\Big).
\end{equation}
\end{prop}
Equation \eqref{prob_bound_1_xi} guarantees that, with probability approaching one, $(\widehat{\xi}_t)$ is bounded away from zero, hence the normalization based on~$(\widehat{\xi}_t)$ is well-behaved. Furthermore, equation~\eqref{sum_bound_1_xi} guarantees the consistency of~$(\widehat{\xi}_t)$ in estimating $({\xi}_t)$.  We need both $p, n\to\infty$ to estimate $(\xi_t)$ consistently, and $p$ can be much larger than $n$.

The next theorem gives the convergence rate of the IPSN pilot estimators, $\widehat{\boldsymbol{\Sigma}}_0$, $\widehat{\boldsymbol{\Lambda}}_{0K}$, and~$\widehat{\boldsymbol{\Gamma}}_{K}$.
\begin{thm}\label{prop:pilot}Under Assumptions \ref{asump1}--\ref{assump_factor_structure}, if in addition, $p, n\to \infty$ and satisfy \\$(\log(p))^{2+\gamma}=o(n)$ for some $\gamma>0$ and $\log(n)=O(\log (p))$, then
\begin{equation}\label{pilot:Sigma_0_all}
\aligned
&\|\widehat{\boldsymbol{\Sigma}}_0-\boldsymbol{\Sigma}_0\|_{\max}=O_p\Bigg(\sqrt{\frac{\log p}{p}}+\sqrt{\frac{\log p}{n}}\Bigg),\\
&\|\widehat{\boldsymbol{\Lambda}}_{0K}\boldsymbol{\Lambda}_{0K}^{-1}-\mathbf{I}\|_{\max}=O_p\Bigg(\sqrt{\frac{\log p}{p}}+\sqrt{\frac{\log p}{n}}\Bigg),\quad \text{and}\\
&\|\widehat{\boldsymbol{\Gamma}}_{K}-\boldsymbol{\Gamma}_{K}\|_{\max}=O_p\Bigg(\sqrt{\frac{\log p}{p^2}}+\sqrt{\frac{\log p}{np}}\Bigg).
\endaligned
\end{equation}
\end{thm}
Theorem \ref{prop:pilot} states that the pilot estimators of the scatter matrix achieve the  sub-Gaussian convergence rates under $2+\varepsilon$th moment assumption.
By contrast, when only $2+\varepsilon$th  moment exists, the pilot estimator $\widehat{\boldsymbol{\Sigma}}$ in \cite{fan2018large} has  a convergence rate of only $O_p(((\log p)/n)^{\varepsilon/(2+\varepsilon)})$.

\subsubsection{Convergence Rates of Generic POET Estimators based on IPSN}
Under the factor model specified in Assumption \ref{assump_factor_structure}, we can write $\boldsymbol{\Sigma}_0=\mathbf{B}_0 \mathbf{B}_0^T+\boldsymbol{\Sigma}_{0u}$, where $\mathbf{B}_0$ is a $p\times K$ matrix, and  $\boldsymbol{\Sigma}_{0u}$ is the idiosyncratic component of the scatter matrix $\boldsymbol{\Sigma}_0$.

We consider factor models with conditional sparsity.  Specifically, we assume that~$\boldsymbol{\Sigma}_{0u}$ belongs to the following  class of sparse matrices: for some  $q\in [0,1]$,  $c>0$, and~\mbox{$s_0(p)<\infty$},
\begin{equation}\label{def:sparse}
\aligned
&\mathcal{U}_q\big(s_0(p)\big)\\
=&\Big\{\boldsymbol{\Sigma}:\boldsymbol{\Sigma} \text{ is positive semi-definite},\, \max_{i}\sigma_{ii}\leq c, \text{and } \max_i\sum_{j=1}\sigma_{ij}^q\leq s_0(p)\Big\}.
\endaligned
\end{equation}
\begin{assumps}\label{sparse_u}
$\boldsymbol{\Sigma}_{0u}\in \mathcal{U}_q\big(s_0(p)\big)$ for some $q\in [0,1]$.
\end{assumps}
Denote $\omega_n=\sqrt{(\log p)/p}+\sqrt{(\log p)/n}$. The next theorem gives the convergence rate of our proposed generic POET estimator of the scatter matrix.
\begin{thm}\label{thm2}
Under the assumptions of Theorem \ref{prop:pilot} and Assumption \ref{sparse_u}, if  in addition, $s_0(p)\omega_n^{1-q}=o(1)$, then
\begin{equation}\label{gpoet_sigma0}
\aligned
&\|(\widehat{\boldsymbol{\Sigma}}_0^{\tau})^{-1}-(\boldsymbol{\Sigma}_0)^{-1}\|_{2}=O_p\Big(s_0(p)\omega_n^{1-q}\Big),\\
&\|\widehat{\boldsymbol{\Sigma}}_0^{\tau}-\boldsymbol{\Sigma}_0\|_{\boldsymbol{\Sigma}_0}=O_p\big(\sqrt{p}\omega_n^2+s_0(p)\omega_n^{1-q}\Big),\quad \text{and}\\
&\|\widehat{\boldsymbol{\Sigma}}_{0 u}^{\tau}-\boldsymbol{\Sigma}_{0u}\|_2=O_p\Big(s_0(p)\omega_n^{1-q}\Big)=\|(\widehat{\boldsymbol{\Sigma}}_{0 u}^{\tau})^{-1}-\boldsymbol{\Sigma}_{0u}^{-1}\|_{2}.
\endaligned
\end{equation}
\end{thm}
Theorem \ref{thm2} asserts that in estimating of the scatter matrix, its idiosyncratic components and their inverse matrices,  under only $2+\varepsilon$th moment condition, our generic~POET estimator achieves the same convergence rate as the sub-Gaussian case; see Theorems 3.1 and 3.2  in \cite{fan2013large}.

Finally, about the estimation of the covariance matrix, recall that we combine the estimators $\widehat{E(\xi_t^2)}$ and $\widehat{\boldsymbol{\Sigma}}_0^{\tau}$ to estimate $\boldsymbol{\Sigma}$, and combine $\widehat{E(\xi_t^2)}$ and $\widehat{\boldsymbol{\Sigma}}_{0u}^{\tau}$ to estimate $\boldsymbol{\Sigma}_u$.  The next proposition gives the convergence rate of $\widehat{E(\xi_t^2)}$ defined in equation \eqref{xi2_est} in estimating~$E(\xi_t^2)$.
\begin{prop}\label{Consist:v} Under Assumptions \ref{asump1}--\ref{assump_factor_structure}, if in addition, $n\to \infty$ and $\log (p)=O(n)$, then
$$
|\widehat{E(\xi_t^2)}-E(\xi_t^2)|=O_p\Bigg(\frac{1}{n^{\varepsilon/(2+\varepsilon)}}+\sqrt{\frac{\log p}{n}}\Bigg).
$$
\end{prop}
The error term $1/n^{\varepsilon/(2+\varepsilon)}$ is from robust Huber estimation based on the infeasible series $(\xi_t)$, and the error term $\sqrt{\log p/n}$ comes from the estimation error of $\widehat{\boldsymbol{\mu}}$. When $\log p=o(n^{(2-\varepsilon)/(2+\varepsilon)})$, the convergence rate of $\widehat{E(\xi_t^2)}$ is the same as the robust Huber estimator applied to  the unobserved series $(\xi_t)$.

The next theorem gives the convergence rates of  $\widehat{\boldsymbol{\Sigma}}^{\tau}$ and $\widehat{\boldsymbol{\Sigma}}_{u}^{\tau}$ in estimating the covariance matrix  and the idiosyncratic covariance matrix, respectively.

\begin{thm}\label{thm3}
Under the assumptions of Theorem \ref{prop:pilot} and Assumption \ref{sparse_u}, if in addition, $s_0(p)\omega_n^{1-q}=o(1)$, then
\begin{equation}\label{gpoet_sigma1}
\aligned
&\|(\widehat{\boldsymbol{\Sigma}}^{\tau})^{-1}-\boldsymbol{\Sigma}^{-1}\|_{2}=O_p\Bigg(s_0(p)\omega_n^{1-q}+\frac{1}{n^{\varepsilon/(2+\varepsilon)}}\Bigg),\\
&\|\widehat{\boldsymbol{\Sigma}}^{\tau}-\boldsymbol{\Sigma}\|_{\boldsymbol{\Sigma}}=O_p\Bigg(\sqrt{p}\omega_n^2+s_0(p)\omega_n^{1-q}+\frac{1}{n^{\varepsilon/(2+\varepsilon)}}\Bigg),\quad \text{and}\\
&\|\widehat{\boldsymbol{\Sigma}}_{u}^{\tau}-\boldsymbol{\Sigma}_{u}\|_2=O_p\Bigg(s_0(p)\omega_n^{1-q}+\frac{1}{n^{\varepsilon/(2+\varepsilon)}}\Bigg)=\|(\widehat{\boldsymbol{\Sigma}}_{u}^{\tau})^{-1}-\boldsymbol{\Sigma}_{u}^{-1}\|_{2}.
\endaligned
\end{equation}
\end{thm}
One can show that when only  $2+\varepsilon$th moment exists, the generic POET estimator proposed in \cite{fan2018large} has a lower convergence rate. Specifically,  it has a convergence rate with the term $\omega^2_n$ replaced by $\widetilde{\omega}^2_n$ and the term $\omega_n^{1-q}$ replaced by $\widetilde{\omega}_n^{2\varepsilon(1-q)/(2+\varepsilon)}$ for $\widetilde{\omega}_n=\sqrt{\log p/n}+1/\sqrt{p}$.

\section{Simulation Studies}\label{Simu}
\subsection{Simulation Setting}
We generate data from the following model. We first simulate~$\mathbf{B}=(b_{ik})_{1\leq i\leq N, 1\leq k\leq K}$, where $b_{ik}\underset{\text{i.i.d.}}\sim N(0, s_i)$,  $K=3$, $s_1=1$, $s_2=0.75^2$, and~\mbox{$s_3=0.5^2$}. We then set $\boldsymbol{\Sigma}_0=p(\mathbf{B}\mathbf{B}^T+\mathbf{I}_p) /\tr ((\mathbf{B}\mathbf{B}^T+\mathbf{I}_p))$.

 We generate the scalar process $(\xi_t)$ from  a Pareto distribution with $P(\xi_t>x)=(x_m/x)^\alpha$ for $x>x_m$, where $x_m$ is chosen such that $E(\xi_t^2)=1$. We set $\alpha=2+\varepsilon$ with $\varepsilon=2$ or $0.2$.  When $\varepsilon=2$, $(\xi_t)$ has finite moments with order below~4, and when $\varepsilon=0.2$, $(\xi_t)$ does not have a finite third moment.

We then generate  $(\mathbf{y}_{t})$ from {$\mathbf{y}_{t}=\xi_t\sqrt{p}\boldsymbol{\Sigma}_0^{1/2}\mathbf{z}_t/\|\mathbf{z}_t\|_2$}, $t=1, ..., n$, where $\mathbf{z}_t\underset{\text{i.i.d.}}\sim N(0, \mathbf{I})$. {Note that the estimators of the scatter matrices and covariance matrices are location invariant, hence without lost of generality, we set $\boldsymbol{\mu}=\mathbf{0}$.}   The dimension and the sample size $p$ and $n$ are set to be $p=500$ and $n=250$, respectively. 

\subsection{Covariance Matrix Estimators and Evaluation Metrics}
{We evaluate the performance of the following covariance matrix estimators: the generic POET estimators based on our idiosyncratic-projected self-normalization approach, denoted by GPOET-IPSN, the robust generic POET estimators in \cite{fan2018large},  denoted by GPOET-FLW, and the original POET estimator based on the sample covariance matrix (\cite{fan2013large}), denoted by POET-S.

We use the same evaluation metrics as  in \cite{fan2018large}. Specifically, errors in estimating $\boldsymbol{\Sigma}_0$ and $\boldsymbol{\Sigma}$ are measured by the difference in the relative Frobenius norm:$$
\aligned&\|\widehat{\boldsymbol{\Sigma}}^{\tau}_0-{\boldsymbol{\Sigma}}_0\|_{{\boldsymbol{\Sigma}}_0}=\frac{1}{\sqrt{p}} \|{\boldsymbol{\Sigma}}_0^{-1/2}(\widehat{\boldsymbol{\Sigma}}^{\tau}_0-{\boldsymbol{\Sigma}}_0){\boldsymbol{\Sigma}}_0^{-1/2}\|_{F},\\
&\|\widehat{\boldsymbol{\Sigma}}^{\tau}-{\boldsymbol{\Sigma}}\|_{{\boldsymbol{\Sigma}}}=\frac{1}{\sqrt{p}}  \|{\boldsymbol{\Sigma}}^{-1/2}(\widehat{\boldsymbol{\Sigma}}^{\tau}-{\boldsymbol{\Sigma}})\widehat{\boldsymbol{\Sigma}}^{-1/2}\|_{F}.
\endaligned$$
Errors in estimating the  matrices $\boldsymbol{\Sigma}_{0u}$, $\boldsymbol{\Sigma}_u$, $\boldsymbol{\Sigma}_{0}^{-1}$, and $\boldsymbol{\Sigma}^{-1}$  are measured  by the difference in spectral norm. }

\subsection{Simulation Results}\label{simu_result}

\begin{figure}
\centering
\includegraphics[width=0.9\textwidth]{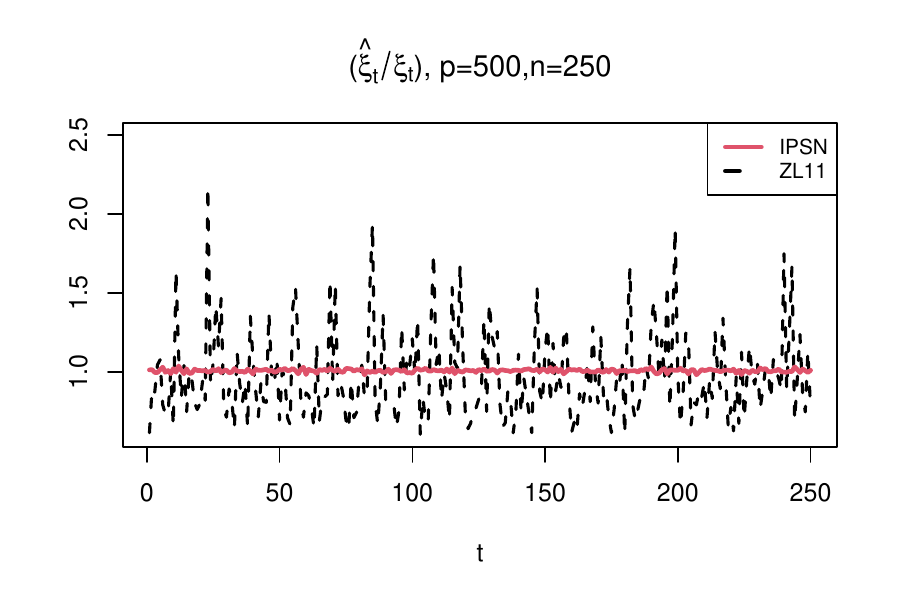}
\caption{{The ratios between $(\widehat{\xi}_t)$ and $(\xi_t)$ from one simulation. We compare our proposed IPSN estimator with the estimator (ZL11) in \cite{ZL11}. The process  $(\xi_t)$ is generated from Pareto distribution $P(\xi_t>x)=(x_m/x)^\alpha$ with the tail parameter $\alpha=2.2$.  }}\label{Fig1}
\end{figure}

First, to illustrate how well our approach works in estimating the scalar process $(\xi_t)$, we compute the estimator $(\widehat{\xi}_t)$ from our  IPSN approach in equation \eqref{hatxi_t} and obtain the ratio $(\widehat{\xi}_t/\xi_t)$ from one realization. We compare the results with the estimator in \cite{ZL11} (denoted by ZL11).
The results are shown in~Figure~\ref{Fig1}.
\noindent We see from Figure~\ref{Fig1} that the ratios between our estimated  and the true $(\xi_t)$'s are remarkably close to one for all $t$. In contrast, due to the factor structure in the data, the estimator  in \cite{ZL11} is not consistent, reflected in that the ratios between the estimated $(\xi_t)$ process and the true ones can be far away from one.


Next, we summarize the performance of the pilot estimators. We include our  IPSN pilot estimators and the following benchmark methods for comparison: the robust estimators of \cite{fan2018large} (FLW), as well as the sample covariance matrix with its leading eigenvectors and eigenvalues (SAMPLE).
For the benchmark methods, we get the pilot estimators for $\mathbf{\Sigma}_0$ and $\mathbf{\Lambda}_{0K}$ by normalizing~$\widehat{\boldsymbol{\Sigma}}$ and $\widehat{\mathbf{\Lambda}}_{K}$ with  $p/\tr(\widehat{\boldsymbol{\Sigma}})$.
  For the pilot estimators of $\widehat{\boldsymbol{\Sigma}}_0$ and $\widehat{\boldsymbol{\Sigma}}$,
we compute the error in the maximum norm.
For the estimators of~$\mathbf{\Gamma}_{K}$, we evaluate the maximum error scaled by $\sqrt{p}$, $(\|\widehat{\mathbf{\Gamma}}_{K}-\mathbf{\Gamma}_{K}\|_{\max})\sqrt{p}.$
For the estimators of $\mathbf{\Lambda}_{0,K}$ and $\mathbf{\Lambda}_{K}$, we report the errors in ratios: $\|\widehat{\mathbf{\Lambda}}_{0,K}\mathbf{\Lambda}_{0K}^{-1}-\mathbf{I}\|_{\max}, \quad\text{and}\quad \|\widehat{\mathbf{\Lambda}}_{K}\mathbf{\Lambda}_{K}^{-1}-\mathbf{I}\|_{\max}.$ In Table \ref{simu_est_pilot_esp1}, we report the  performance of the pilot estimators from 100 replications.
\begin{table}[H]
\caption{Performance of the pilot estimators. The heavy-tailedness index is set to be $\varepsilon=2$ or $0.2$. The dimension and the sample size are $p=500$ and $n=250$, respectively. Reported values are mean and standard deviations (in parentheses) from~100 replications for the error in maximum norm in estimating $\boldsymbol{\Sigma}_0$, $\boldsymbol{\Sigma}$, and ($\sqrt{p}\times$) $\boldsymbol{\Gamma}_{K}$, and the maximum norm error in ratios in estimating $\boldsymbol{\Lambda}_{0K}$ and $\boldsymbol{\Lambda}_{K}$. We include our IPSN and the following benchmark methods: the robust method proposed in \cite{fan2018large} (FLW) and the sample covariance matrix (SAMPLE).
}
\begin{center}
\tabcolsep 0.1in\renewcommand{\arraystretch}{1.2} \doublerulesep
2pt
\begin{tabular}{lcccccccc}
\hline \hline
&$\boldsymbol{\Sigma}_0$&$\boldsymbol{\Sigma}$&$\boldsymbol{\Lambda}_{0K}$&$\boldsymbol{\Lambda}_{K}$&$\boldsymbol{\Gamma}_{K}$\\\hline
$\varepsilon=2$\\
IPSN&	0.482	&	0.676	&	0.100	&	0.165	&	0.592		\\
&(0.108)&(0.273)&(0.048)&(0.099)&(0.177)\\
FLW&	0.706	&	0.921	&	0.111	&	0.181	&	0.620		\\
&(0.392)&(0.760)&(0.053)&(0.122)&(0.178)\\
SAMPLE&	0.728	&	0.923	&	0.135	&	0.204	&	0.872		\\
&(0.394)&(0.768)&(0.076)&(0.147)&(0.488)\\\\
$\varepsilon=0.2$\\
IPSN&0.485	&	2.105	&	0.100	&	0.623	&	0.598	\\
&(0.109)&(2.706)&(0.048)&(0.886)&(0.182)\\
FLW&1.820	&	3.371	&	0.141	&	0.656	& 0.641	\\
&(1.289)&(9.982)&(0.062)&(0.891)&(0.185)\\
SAMPLE&1.913	&	3.406&	0.326	&	0.882	&	2.135\\
&(1.313)&(9.987)&(0.248)&(2.350)&(1.147)\\
\hline
\end{tabular}
\end{center}\label{simu_est_pilot_esp1}
\end{table}
\noindent We see from Table \ref{simu_est_pilot_esp1} that our  IPSN estimators outperform the other estimators in all measures.
 The advantage of IPSN over~FLW and SAMPLE is  more salient in the estimation of the scatter matrices and with increasing heavy-tailedness in the data. In particular,
  for the pilot estimators of the  scatter matrix and its leading eigenvectors and leading eigenvalues, when~$\varepsilon$ changes from~2 to 0.2, that is, when data get more  heavy-tailed, the estimation errors for~IPSN remain almost the same. By contrast, the errors in the FLW and SAMPLE estimators increase substantially.
For pilot estimators of the covariance matrix, the error of our IPSN estimator grows with the heavy-tailedness of the model, but it still significantly outperforms the other methods.


Finally, we summarize the  performance of the generic POET estimators in  Table~\ref{simu_est_gpoet_sig_esp1}.
\begin{table}[H]
\caption{Performance of generic POET estimators. The heavy-tailedness index is set to be $\varepsilon=2$ or $0.2$. The dimension and the sample size are $p=500$ and $n=250$, respectively. Reported values are
mean and standardized deviation (in parentheses) from~100 replications  for the error in relative Frobenius norm in estimating $\boldsymbol{\Sigma}_0$  and~$\boldsymbol{\Sigma}$, and the error in $\ell_2$ norm in estimating $\boldsymbol{\Sigma}_{0u}$,  $\boldsymbol{\Sigma}_{0u}$, $\boldsymbol{\Sigma}_{0}^{-1}$ and $\boldsymbol{\Sigma}^{-1}$. We include  our GPOET-IPSN and the following benchmark methods: the robust method GPOET-FLW and the original POET based on sample covariance matrix (POET-S).
}
\begin{center}
\tabcolsep 0.05in\renewcommand{\arraystretch}{1.2} \doublerulesep
2pt
\begin{tabular}{lcccccccccc}
\hline \hline
&$\boldsymbol{\Sigma}_0$&$\boldsymbol{\Sigma}$&$\boldsymbol{\Sigma}_{0u}$&$\boldsymbol{\Sigma}_{u}$&$\boldsymbol{\Sigma}_0^{-1}$&$\boldsymbol{\Sigma}^{-1}$&$\boldsymbol{\Sigma}_{0u}^{-1}$&$\boldsymbol{\Sigma}_{u}^{-1}$\\\hline
$\varepsilon=2$\\
GPOET-IPSN&0.241&0.091&0.107&0.122&1.023&2.894&1.029&1.161\\
&(0.008)&(0.020)&(0.012)&(0.034)&(0.209)&(0.396)&(0.210)&(0.301)\\
GPOET-FLW&0.328&0.127&0.545&0.632&1.533&2.995&1.542&1.621\\
&(0.103)&(0.085)&(0.768)&(1.360)&(0.354)&(0.290)&(0.357)&(0.318)\\
POET-S&0.479&0.186&0.206&0.224&1.423&3.235&1.430&1.494\\
&(0.483)&(0.278)&(0.239)&(0.318)&(0.437)&(0.284)&(0.439)&(0.289)\\\\
$\varepsilon=0.2$\\
GPOET-IPSN&0.242&0.219&0.107&0.266&1.023&7.155&1.028&5.448\\
&(0.009)&(0.289)&(0.013)&(0.371)&(0.210)&(2.567)&(0.211)&(2.417)\\
GPOET-FLW&0.707&0.375&2.245&3.098&4.822&11.321&6.839&12.611\\
&(0.394)&(1.218)&(2.428)&(15.229)&(10.292)&(10.969)&(20.067)&(27.690)\\
POET-S&3.331&1.381&0.549&0.412&4.925&11.616&4.972&9.888\\
&(4.122)&(5.096)&(0.567)&(0.434)&(11.121)&(2.812)&(11.203)&(2.996)\\\hline
\end{tabular}
\end{center}\label{simu_est_gpoet_sig_esp1}
\end{table}
\noindent We see from Table \ref{simu_est_gpoet_sig_esp1} that
our GPOET-IPSN estimators deliver the lowest estimation error in all measures.
Compared to GPOET-FLW or POET-S, the advantage of GPOET-IPSN is especially evident  for the estimation of $\boldsymbol{\Sigma}_0$, $\boldsymbol{\Sigma}_{0u}$, their inverse matrices, and for the more heavy-tailed setting when $\varepsilon=0.2$. When estimating the scatter matrix, its idiosyncratic component and their inverses, the performance of GPOET-IPSN is robust for different heavy-tailedness conditions. When~$\varepsilon$ decreases from 2 to~0.2, the data become more heavy-tailed, but the estimation errors of GPOET-IPSN remain largely the same.
By contrast,  the errors of the GPOET-FLW and POET-S estimators increase sharply in  the more heavy-tailed setting.
In the estimation of $\boldsymbol{\Sigma}$ and $\boldsymbol{\Sigma}^{-1}$, GPOET-IPSN also has  lower estimation errors than the benchmark methods.


In summary, the simulation results validate the  theoretical properties of our approach and demonstrate its substantial advantage in estimating the scatter matrix and its inverse matrix. Such advantages are particularly  valuable  in  applications where only the scatter matrix matters while  the scalar component does not, as we will show in the empirical studies below.

\section{Empirical Studies}\label{Empi}
\subsection{Global Minimum Variance Portfolio Optimization}
The estimation of the minimum variance portfolio (MVP) in the high-dimensional setting has drawn considerable attention in recent years. It is widely used to measure the empirical performance of covariance matrix estimators (\cite{fan2012vast, ledoit2017nonlinear, DLZ19}). In this subsection, we conduct empirical analysis to evaluate the performance of our proposed estimator in  the MVP optimization.

The MVP is the portfolio that achieves the minimum variance among all portfolios with the constraint that the summation of weights equals one. Specifically, for $p$ asset returns with a covariance matrix of $\boldsymbol{\Sigma}$, the MVP solves
$$
\mathop{\text{min}}\mathbf{w}^T\boldsymbol{\Sigma}\mathbf{w}, \text{ subject to }\mathbf{w}^T\mathbf{1}=1,
$$
where $\mathbf{w}$ is a length $p$ vector of portfolio weights, and $\mathbf{1}$ is a length $p$ vector of ones. The solution of the MVP is
$$
\mathbf{w}^*=\frac{1}{\mathbf{1}^T\boldsymbol{\Sigma}^{-1}\mathbf{1}}\boldsymbol{\Sigma}^{-1}\mathbf{1}.
$$
It is straightforward to see that under elliptical models,
$$
\mathbf{w}^*=\frac{1}{\mathbf{1}^T\boldsymbol{\Sigma}_0^{-1}\mathbf{1}}\boldsymbol{\Sigma}_0^{-1}\mathbf{1}.
$$
In other words,  using ${\boldsymbol{\Sigma}}^{-1}$ and ${\boldsymbol{\Sigma}}_0^{-1}$ are equivalent in constructing the MVP.
Using an estimator of the inverse of the scatter matrix or the covariance matrix, $\widehat{\boldsymbol{\Sigma}}_0^{-1}$ or $\widehat{\boldsymbol{\Sigma}}^{-1}$,  we estimate the~MVP by
\begin{equation}\label{hatw}
\widehat{\mathbf{w}}^*=\frac{1}{\mathbf{1}^T\widehat{\boldsymbol{\Sigma}}_0^{-1}\mathbf{1}}\widehat{\boldsymbol{\Sigma}}_0^{-1}\mathbf{1},\quad \text{or}\quad \widehat{\mathbf{w}}^*=\frac{1}{\mathbf{1}^T\widehat{\boldsymbol{\Sigma}}^{-1}\mathbf{1}}\widehat{\boldsymbol{\Sigma}}^{-1}\mathbf{1}.
\end{equation}

\subsection{Data and Compared Methods}
We use the daily returns of S\&P 500 Index constituent stocks between January 1995 and December 2023 to construct MVPs. At the beginning of each month, we use the historical returns of the stocks that stayed in the S\&P 500 Index for the past five years to estimate the portfolio weights. We evaluate the  risk of the portfolios based on their out-of-sample daily returns from January~2000 to December~2023.

We estimate the MVP using equation \eqref{hatw} with the following estimators:  our  GPOET estimator of the scatter matrix based on idiosyncratic-projected self-normalization, \mbox{GPOET-IPSN}, the GPOET-FLW estimator, and the original POET estimator based on the sample covariance matrix, POET-S. When performing the generic POET procedure, the tuning parameter is chosen by cross-validation with the criterion of minimizing the out-of-sample MVP risk. The number of factors is estimated by the method in \cite{yu2019robust}.  We also include the equal-weight portfolio $(1/p, ..., 1/p)^T$ as a benchmark portfolio, denoted by EW.

\subsection{Out-of-Sample Performance}
We use  the standard deviation of out-of-sample returns to evaluate the out-of-sample risk of portfolios.  In addition, we perform  the following test to evaluate the statistical significance of the differences in risks between our GPOET-IPSN method and benchmark portfolios:
\begin{equation}\label{test}
H_0: \sigma\geq \sigma_{0}\quad\text{ vs. }\quad H_1: \sigma< \sigma_{0},
\end{equation}
where~$\sigma$ denotes the standard deviation of the GPOET-IPSN portfolio and~$\sigma_{0}$ denotes the standard deviation of a benchmark portfolio. We adopt the test in \cite{LW11} and use the heteroskedasticity-autocorrelation-consistent~(HAC) standard deviation estimator of the test statistic therein.

Table \ref{emp_risk} reports the risks of various portfolios.

 \begin{table}[H]
\caption{Out-of-sample risks of the portfolios. We report the annualized standard deviations. We construct the MVP portfolio  using equation \eqref{hatw} with the covariance matrix  estimated using generic POET based on our idiosyncratic-projected self-normalized method (GPOET-IPSN), the robust method GPOET-FLW, and the sample covariance matrix (POET-S). We also include the equal weight (EW) portfolio as a benchmark. The out-of-sample period is from January 2000 to December 2023. The risk is evaluated using out-of-sample daily returns of the portfolios. For all methods, the portfolios weights are estimated based on prior five-year daily excess returns of S\&P 500 Index constituent stocks, reestimated monthly. The symbols~*, ** and *** indicate statistical significance for testing~\eqref{test} at  5\%, 1\% and 0.1\% levels, respectively.
}
\begin{center}
\tabcolsep 0.1in\renewcommand{\arraystretch}{1.1} \doublerulesep
2pt
\begin{tabular}{llllll}
\hline \hline
&\multicolumn{5}{c}{Out-of-sample portfolio risk}\\\hline
Period& 2000--2023 &2000  &2001  &2002&2003\\
GPOET-IPSN &0.121 &0.154 &0.116 &0.140 &0.090\\
GPOET-FLW &0.150 *** &0.168 **&0.167 ***&0.173 ***&0.121 ***\\
POET-S &0.150 *** &0.152&0.121  &0.140&0.095 **\\
{EW}  &	0.217	***	&0.192 ***&0.208 *** &0.269 ***&	0.179 ***\\\\
Period& 2004&2005 &2006&2007 &2008\\
GPOET-IPSN &0.081&0.084&0.068 &0.081&0.208\\
GPOET-FLW &0.095 ***&0.100 *** &0.089 ***&0.100 ***&0.224 \\
POET-S &0.091 ***&0.091 *** &0.074**&0.087 ***&0.188\\
{EW}  &0.128 ***&0.115 *** &0.113 ***&	0.166 ***&0.460 ***\\\\
Period   &2009 &2010 &2011 &2012 &2013\\
 GPOET-IPSN &0.133&0.083&0.103&0.082&0.084\\
GPOET-FLW *&0.228 *** &0.094 **&0.108&0.073 &0.081 \\
POET-S &0.189 *** &0.089 * &0.114 **&0.084 &0.089\\
{EW}  &0.357 *** &0.210 ***&	0.271 ***&0.148***&0.123 ***\\\\
Period &2014 &2015 &2016&2017 &2018\\
 GPOET-IPSN &0.076&0.096&0.099&0.072&0.105\\
GPOET-FLW &0.081 &0.114 ***  &0.128 ***&0.076 &0.116 *\\
POET-S &0.089 &0.120 ***&0.116 **& 0.073 &0.119 ***\\
{EW}   &0.119 ***&	0.158 ***&0.157 ***&0.077 &0.161 ***\\\\
Period  &2019 &2020 &2021 &2022 &2023\\
 GPOET-IPSN &0.097&0.269&0.106&0.141&0.125\\
GPOET-FLW&0.131 *** &0.317 **&0.184 *** &0.179 **&0.176 ***\\
POET-S&0.101 &0.413 ***&0.194 *** &0.226 ***&0.168 ***\\
{EW}  &	0.137 ***&0.404 ***&0.144 *** &0.239 *** &0.152 ***\\\hline
\end{tabular}
\end{center}\label{emp_risk}
\end{table}
\noindent We see from Table \ref{emp_risk} that our GPOET-IPSN portfolio  achieves a substantially lower risk  than benchmark  portfolios for the period between 2000 and 2023. The statistical test results suggest that the differences in the risk between our portfolio and benchmark portfolios are all  statistically significantly negative. In addition, when checking the risk performance in different years, we see that the  GPOET-IPSN portfolio performs robustly well. It outperforms the benchmark portfolios by achieving the lowest risk in the majority of the years (20 out of 24 years). The statistical test results show that the difference between the risks of our GPOET-IPSN  portfolio and other portfolios is statistically significantly negative for most of the years.

\section{Conclusion}\label{Conc}

We propose an idiosyncratic-projected self-normalization (IPSN) approach to estimate the scatter matrix under high-dimensional elliptical factor models. We show that our estimator achieves the sub-Gaussian convergence rate  under only $2+\varepsilon$th moment condition, a property that can not be achieved by existing approaches. Moreover, we develop  POET estimators of the scatter matrix and the covariance matrix under a conditional sparsity assumption, and show that they have higher convergence rates than the generic POET estimator in \cite{fan2018large}.  Numerical studies demonstrate the clear advantages of  our proposed estimators over various benchmark methods.
\section{Proofs}\label{Proof}
In this section we prove Theorems \ref{prop:pilot}--\ref{thm3} and Propositions \ref{lem_denominator_xi}--\ref{Consist:v}. We first give the following lemmas. The proofs of Lemmas \ref{mu_est}--\ref{prop:pilot:Sigma} are given in the Supplementary Material (\cite{DZ24_supp}).

\subsection{Lemmas}
\begin{lem}\label{huber}[Corollary 1 of \cite{avella2018robust}] Suppose that $\delta\in (0, 1)$, $\varepsilon \in(0, 1]$,  $n>12\log(2\delta^{-1})$, and  $Z_1$, ..., $Z_n$ are \mbox{i.i.d.} random variables with mean $g$ and bounded $1+\varepsilon$th moment, i.e. $E(|Z_1-g|^{1+\varepsilon})=v<\infty$. Take $H=\big(vn/\log(2\delta^{-1})\big)^{1/(1+\varepsilon)}$. Let $\widehat{g}_H$ satisfies that $\sum_{t=1}^n\Psi_H(Z_t-\widehat{g}_H)=0$ with $\Psi_H$ being the Huber function defined in equation \eqref{Psi_H}. Then with probability at least $1-2\delta$,
$$
|\widehat{g}_H-g|\leq C v^{1/(1+\varepsilon)}\Big(\frac{\log (2\delta^{-1})}{n}\Big)^{\varepsilon/(1+\varepsilon)}.$$
\end{lem}

\begin{lem}\label{ZL2012}[Lemma 9 of \cite{ZL11}]  Suppose that $\mathbf{X}=(x_1, ..., x_p)^T$, where $x_i$'s are \mbox{i.i.d.} random variables such that $E(x_1)=0$, $E(x_1^2)=1$ and $E(x_1^{2k})<\infty$ for some $2\leq k<\infty$. Then there exists $c_k>0$, depending only on $k$, $E(x_1^4)$ and $E(x_1^{2k})$ such that for any $p\times p$ nonrandom matrix $\mathbf{A}$,
$$
E\Big(|\mathbf{X}^T\mathbf{A}\mathbf{X}-\tr (\mathbf{A})|^{2k}\Big)\leq c_k \Big(\tr(\mathbf{A}\mathbf{A}^T)\Big)^k\leq c_kp^k\|\mathbf{A}\|_2^{2k}.
$$
\end{lem}
\begin{lem}\label{Chatterjee2013}[Lemma 3 of \cite{chatterjee2013assumptionless}]  Suppose that $x_i\sim N(0, \sigma_i^2)$, $i=1, ..., m$. The $x_i$'s need not be independent. Let $L:=\max_{1\leq i\leq m}\sigma_i$. Then
$$
E\Big(\max_{1\leq i\leq m}|x_i|\Big)\leq L\sqrt{2\log(2m)}.
$$
\end{lem}

\begin{lem}\label{mu_est}Under Assumptions \ref{asump1} and \ref{sigma_bound}, if in addition $\log(p)=O(n)$, then
$$
\|\widehat{\boldsymbol{\mu}}-\boldsymbol{\mu}\|_{\max}=O_p\Bigg(\sqrt{\frac{\log p}{n}}\Bigg).$$
\end{lem}


\begin{lem}\label{normal_quadraticratios} Suppose that $\mathbf{z}_t=(z_1, ..., z_p)^T\sim {N}(0, \mathbf{I})$. For any $p\times p$ positive semi-definite matrix $\mathbf{\Sigma}$, we have
\begin{equation}\label{E_ratio_trace}
E\Big(\frac{\mathbf{z}_t^T\boldsymbol{\Sigma}\mathbf{z}_t}{\mathbf{z}^T_t\mathbf{z}_t}\Big)=\frac{\tr(\boldsymbol{\Sigma})}{p},
\end{equation}
and
\begin{equation}\label{E_ratio2_trace}
E\Bigg(\frac{(\mathbf{z}_t^T\boldsymbol{\Sigma}\mathbf{z}_t)^2}{(\mathbf{z}^T_t\mathbf{z}_t)^2}\Bigg)\leq \frac{\tr(\boldsymbol{\Sigma}^2)}{p}+\frac{\big(\tr(\boldsymbol{\Sigma})\big)^2}{p(p-1)}.
\end{equation}
\end{lem}

Define
$$\widetilde{\mathbf{X}}_t=\boldsymbol{\Sigma}_0^{1/2}\mathbf{z}_t, \quad \text{ and }\quad \widetilde{\mathbf{S}}_0=\frac{1}{n}\sum_{t=1}^n\widetilde{\mathbf{X}}_t\widetilde{\mathbf{X}}_t^T.
$$
Let $\widetilde{\boldsymbol{\Lambda}}_{0K}=\diag(\widetilde{\lambda}_{0;1}, ..., \widetilde{\lambda}_{0;K})$ with $\widetilde{\lambda}_{0;1}\geq ...\geq \widetilde{\lambda}_{0;K}$ be the leading eigenvalues and $\widetilde{\boldsymbol{\Gamma}}_{K}=(\widehat{\phi}_1, ..., \widehat{\phi}_K)$ the corresponding eigenvectors of $\widetilde{\mathbf{S}}_0$.
\begin{lem}\label{gaussian_pilot}Under Assumptions \ref{asump1}--\ref{assump_factor_structure}, if in addition,  $\log(p)=o(n)$, then
\begin{equation}\label{Gaussian:Sig}
\|\widetilde{\mathbf{S}}_0-\boldsymbol{\Sigma}_0\|_{\max}=O_p\Bigg(\sqrt{\frac{\log p}{n}}\Bigg),
\end{equation}
\begin{equation}\label{Gaussian:Lam0}
\|\widetilde{\boldsymbol{\Lambda}}_{0K}(\boldsymbol{\Lambda}_{0K})^{-1}-\mathbf{I}\|_{\max}=O_p\Bigg(\sqrt{\frac{\log p}{n}}\Bigg),\quad \text{and}
\end{equation}
\begin{equation}\label{Gaussian:Gam0}
\|\widetilde{\boldsymbol{\Gamma}}_{K}-\boldsymbol{\Gamma}_{K}\|_{\max}=O_p\Bigg(\sqrt{\frac{\log p}{pn}}\Bigg).
\end{equation}
\end{lem}

Recall that $\lambda_{0;i}$ stands for the $i$th largest eigenvalue of $\boldsymbol{\Sigma}_0$. Define
\begin{equation}\label{def:siga}
\sigma_a=\frac{1}{p}\sum_{i=K+1}^p \lambda_{0;i}.
\end{equation}
Under Assumption \ref{assump_factor_structure}, we have that $\sigma_a\asymp 1$.

\begin{lem}\label{lem_denominator}Under Assumptions \ref{asump1}--\ref{assump_factor_structure}, if in addition,  $(\log(p))^{2+\gamma}=o(n)$ for some $\gamma>0$ and $\log(n)=O(\log (p))$, then
\begin{equation}\label{prob_bound_1}
P\Bigg(\min_{1\leq t\leq n}\frac{\|\widehat{\mathbf{P}}_I(\mathbf{y}_t-\widehat{\boldsymbol{\mu}})\|_2^2}{p\xi_t^2}>\frac{\sigma_a}{3}\Bigg)\to 1,
\end{equation}
\begin{equation}\label{sum_bound_1}
\frac{1}{n}\sum_{t=1}^n\Bigg(\frac{\xi_t^2p}{\|\widehat{\mathbf{P}}_I(\mathbf{y}_t-\widehat{\boldsymbol{\mu}})\|_2^2}-\frac{1}{\sigma_a}\Bigg)^2=O_p\Bigg(\frac{\log p}{p}+\frac{\log p}{n}\Bigg),
\end{equation}
and
\begin{equation}\label{trace_x}\aligned
\frac{1}{p}\tr\Big(\frac{1}{n}\sum_{t=1}^n\widehat{\mathbf{X}}_t\widehat{\mathbf{X}}_t^T\Big)=\frac{1}{\sigma_a}+O_p\Bigg(\sqrt{\frac{\log p}{n}}+\sqrt{\frac{\log p}{p}}\Bigg).
\endaligned
\end{equation}
\end{lem}

\begin{lem}\label{prop:pilot:Sigma}Under Assumptions \ref{asump1}--\ref{assump_factor_structure}, if in addition,  $(\log(p))^{2+\gamma}=o(n)$ for some $\gamma>0$ and $\log(n)=O(\log (p))$, then
\begin{equation}\label{pilot:Sigma_0:Sig}
\|\widehat{\boldsymbol{\Sigma}}_0-\widetilde{\mathbf{S}}_0\|_{\max}=O_p\Bigg(\sqrt{\frac{\log p}{p}}+\sqrt{\frac{\log p}{n}}\Bigg),
\end{equation}
\begin{equation}\label{pilot:Sigma_0:Lam0}
\|\widehat{\boldsymbol{\Lambda}}_{0K}(\widetilde{\boldsymbol{\Lambda}}_{0K})^{-1}-\mathbf{I}\|_{\max}=O_p\Bigg(\sqrt{\frac{\log p}{p}}+\sqrt{\frac{\log p}{n}}\Bigg),\quad \text{and}
\end{equation}
\begin{equation}\label{pilot:Sigma_0:Gam0}
\|\widehat{\boldsymbol{\Gamma}}_{K}-\widetilde{\boldsymbol{\Gamma}}_{K}\|_{\max}=O_p\Bigg(\sqrt{\frac{\log p}{p^{2}}}+\sqrt{\frac{\log p}{pn}}\Bigg).
\end{equation}
\end{lem}

\subsection{Proof of Proposition \ref{lem_denominator_xi}} The desired results follow from Lemma \ref{lem_denominator}. $\qed$

\subsection{Proof of Theorem  \ref{prop:pilot}}

The desired bounds follow from Lemmas \ref{gaussian_pilot} and \ref{prop:pilot:Sigma}. $\qed$

\subsection{Proof of Proposition \ref{Consist:v}}
Using the definition of $\widehat{E(\xi_t^2)}$ and the inequality that for any $x\geq y$,
\begin{equation}\label{huber_xy_ineq}
\Psi_H(x)-\Psi_H(y)\geq (x-y)\mathbf{1}_{\big\{x<H \,\text{and}\, y>-H\big\}},
\end{equation} we have
$$
|\widehat{E(\xi_t^2)}-E(\xi_t^2)|\leq \frac{1}{\tilde{p}^H}\Bigg|\frac{1}{n}\sum_{t=1}^n\Psi_H\Big(\|\mathbf{y}_t-\widehat{\boldsymbol{\mu}}\|_2^2/p-E(\xi_t^2)\Big)\Bigg|,
$$
where  $\tilde{p}^H=\sum_{t=1}^n \gamma_t/n$ for $(\gamma_t)$ defined as follows:
\begin{equation*}
\gamma_t=
\begin{cases}
&\mathbf{1}_{\{\|\mathbf{y}_t-\widehat{\boldsymbol{\mu}}\|_2^2/p-\widehat{E(\xi_t^2)}  <H \quad\text{and}\quad\|\mathbf{y}_t-\widehat{\boldsymbol{\mu}}\|_2^2/p-E(\xi_t^2)>-H\}}\quad \text{ if }\widehat{E(\xi_t^2)}\leq E(\xi_t^2)\\
&\mathbf{1}_{\{\|\mathbf{y}_t-\widehat{\boldsymbol{\mu}}\|_2^2/p-\widehat{E(\xi_t^2)}>-H \quad\text{and}\quad\|\mathbf{y}_t-\widehat{\boldsymbol{\mu}}\|_2^2/p-E(\xi_t^2)<H\}}\quad \text{ if }\widehat{E(\xi_t^2)}\geq E(\xi_t^2)
\end{cases}.
 \end{equation*}
We will show that
\begin{equation}\label{Psi_Hbound}
\Bigg|\frac{1}{n}\sum_{t=1}^n\Psi_H\Big(\|\mathbf{y}_t-\widehat{\boldsymbol{\mu}}\|_2^2/p-E(\xi_t^2)\Big)\Bigg|=O_p\Bigg(\sqrt{\frac{\log p}{n}}+\frac{1}{n^{\varepsilon/(2+\varepsilon)}}\Bigg),
\end{equation}
and
\begin{equation}\label{Prob_Hbound}
\tilde{p}^H=1+o_p(1).
\end{equation}
The desired bound will follow from  \eqref{Psi_Hbound} and \eqref{Prob_Hbound}.

We start with \eqref{Psi_Hbound}. By \eqref{Y:model}, we have
\begin{equation}\label{I_III}
\aligned
\frac{\|\mathbf{y}_t-\widehat{\boldsymbol{\mu}}\|_2^2}{p}=&\xi_t^2\frac{\mathbf{z}_t^T\boldsymbol{\Sigma}_0\mathbf{z}_t}{\mathbf{z}^T_t\mathbf{z}_t}+\frac{\|\widehat{\boldsymbol{\mu}}-\boldsymbol{\mu}\|_2^2}{p}+\frac{2}{\sqrt{p}}\xi_t\frac{(\boldsymbol{\mu}-\widehat{\boldsymbol{\mu}})^T\boldsymbol{\Sigma}_0^{1/2}\mathbf{z}_t}{\|\mathbf{z}_t\|_2}\\
=:&I_t+II+III_t.
\endaligned
\end{equation}
The definition of $\Psi_H(\cdot)$ implies that for any $a, b$,
$$
|\Psi_H(a)-\Psi_H(b)|\leq |a-b|.
$$
Hence, for any $(a_i)_{1\leq i\leq n}$, $(b_i)_{1\leq i\leq n}$,
\begin{equation}\label{triangle_ineq}
\Bigg|\sum_{t=1}^n\Psi_H(a_i)-\sum_{t=1}^n\Psi_H(b_i)\Bigg|\leq \sum_{t=1}^n|a_i-b_i|.
\end{equation}
By \eqref{I_III} and \eqref{triangle_ineq}, we get
\begin{equation}\label{partition}
\aligned
&\Bigg|\frac{1}{n}\sum_{t=1}^n\Psi_H\Bigg(\frac{\|\mathbf{y}_t-\widehat{\boldsymbol{\mu}}\|_2^2}{p}-E(\xi_t^2)\Bigg)\Bigg|\\
\leq & \Bigg|\frac{1}{n}\sum_{t=1}^n\Psi_H\Big(I_t-E(\xi_t^2)\Big)\Bigg|+II+\frac{1}{n}\sum_{t=1}^n|III_t|.
\endaligned
\end{equation}
By Lemma \ref{normal_quadraticratios}, and the assumptions that $E(\xi_t^{2+\varepsilon})<\infty$ and   $(\xi_t)$ is independent with $(\mathbf{z}_t)$, we have
\begin{equation}\label{E_I}
E(I_t)=E(\xi_t^2)\text{, and }E(I_t^{1+\varepsilon/2})<\infty.
\end{equation}
By \eqref{E_I} and equation (S31)  in \cite{avella2018robust}, for $H\asymp n^{1/(1+\varepsilon/2)}$,
\begin{equation}\label{part_I}
\Bigg|\frac{1}{n}\sum_{t=1}^n\Psi_H\Big(I_t-E(\xi_t^2)\Big)\Bigg|=O_p\Bigg(\frac{1}{n^{\varepsilon/(2+\varepsilon)}}\Bigg).
\end{equation}
By Lemma \ref{mu_est},
\begin{equation}\label{II}
II=O_p\Bigg(\frac{\log p}{n}\Bigg).
\end{equation}
By the Cauchy-Schwarz inequality,
\begin{equation}\label{III}
\aligned
\frac{1}{n}\sum_{t=1}^n|III_t|\leq& \frac{1}{\sqrt{p}}\sqrt{\frac{1}{n}\sum_{t=1}^n\xi^2_t}\sqrt{\frac{1}{n}\sum_{t=1}^n\frac{\big((\boldsymbol{\mu}-\widehat{\boldsymbol{\mu}})^T\boldsymbol{\Sigma}_0^{1/2}\mathbf{z}_t\big)^2}{\|\mathbf{z}_t\|_2^2}}\\
\leq&\frac{\|\boldsymbol{\mu}-\widehat{\boldsymbol{\mu}}\|_2}{\sqrt{p}}\sqrt{\frac{1}{n}\sum_{t=1}^n\xi^2_t}\sqrt{\frac{1}{n}\sum_{t=1}^n\frac{\mathbf{z}_t^T\boldsymbol{\Sigma}_0\mathbf{z}_t}{\|\mathbf{z}_t\|^2_2}}.
\endaligned\end{equation}
By the assumption that $E(\xi_t^{2+\varepsilon})<\infty$, \eqref{E_ratio_trace} and Markov's inequality,
\begin{equation}\label{part_III1}
\frac{1}{n}\sum_{t=1}^n\xi^2_t=O_p(1),\text{  and }\,\frac{1}{n}\sum_{t=1}^n\frac{\mathbf{z}_t^T\boldsymbol{\Sigma}_0\mathbf{z}_t}{\|\mathbf{z}_t\|^2_2}=O_p(1).
\end{equation}Combining \eqref{II}, \eqref{III} and \eqref{part_III1} yields that
\begin{equation}\label{part_III}
\frac{1}{n}\sum_{t=1}^n|III_t|=O_p\Bigg(\sqrt{\frac{\log p}{n}}\Bigg).
\end{equation}
The desired bound \eqref{Psi_Hbound} follows from \eqref{partition}, \eqref{part_I}, \eqref{II} and \eqref{part_III}.

{Finally, about \eqref{Prob_Hbound}, if $\widehat{E(\xi_t^2)}\leq E(\xi_t^2)$,  we have that for all $n$ large enough,
$$
\aligned
\gamma_t\geq& \mathbf{1}_{\big\{\big|\frac{\|\mathbf{y}_t-\widehat{\boldsymbol{\mu}}\|_2^2}{p}-E(\xi_t^2)\big| \leq \frac{H}{2}\quad\text{and}\quad\widehat{E(\xi_t^2)}\geq 0\big\}}\\
\geq & 1-\mathbf{1}_{\big\{\big|\frac{\|\mathbf{y}_t-\widehat{\boldsymbol{\mu}}\|_2^2}{p}-E(\xi_t^2)\big| >\frac{H}{2}\big\}}-\mathbf{1}_{\big\{\widehat{E(\xi_t^2)} <0 \big\}},
\endaligned
$$
where in the first inequality we use $H\asymp n^{1/(1+\varepsilon/2)}\to\infty$ so that $E(\xi_t^2)<H/2$ for all large~$n$.
Similarly, if $\widehat{E(\xi_t^2)}\geq E(\xi_t^2)$,
$$
\gamma_t\geq 1-\mathbf{1}_{\big\{\big|\frac{\|\mathbf{y}_t-\widehat{\boldsymbol{\mu}}\|_2^2}{p}-E(\xi_t^2)\big| >\frac{H}{2}\big\}}-\mathbf{1}_{\big\{\widehat{E(\xi_t^2)} >\frac{H}{2} \big\}}.
$$
It follows that
\begin{equation}\label{Egammat}
E(\gamma_t)\geq 1-P\Big(\Big|\frac{\|\mathbf{y}_t-\widehat{\boldsymbol{\mu}}\|_2^2}{p}-E(\xi_t^2)\Big|>\frac{H}{2}\Big)-P\Big(\widehat{E(\xi_t^2)}>\frac{H}{2}\Big)-P\Big(\widehat{E(\xi_t^2)}<0\Big).
\end{equation}
We have
$$\aligned
&P\Bigg(\Big|\frac{\|\mathbf{y}_t-\widehat{\boldsymbol{\mu}}\|_2^2}{p}-E(\xi_t^2)\Big|>\frac{H}{2}\Bigg)\\
\leq & P\Bigg(|I_t-E(\xi_t^2)|>\frac{H}{6}\Bigg)+P\Bigg(|II|>\frac{H}{6}\Bigg)+P\Bigg(|III_t|>\frac{H}{6}\Bigg),\endaligned$$
where terms $I_t$, $II$ and $III_t$ are defined in equation  \eqref{I_III}.
By \eqref{E_I}, $H\asymp n^{1/(1+\varepsilon/2)}$,
and Markov' inequality, we have
$$P\Bigg(|I_t-E(\xi_t^2)|>\frac{H}{6}\Bigg)=O\Bigg(\frac{1}{H^{1+\varepsilon/2}}\Bigg)=O\Bigg(\frac{1}{n}\Bigg)\to 0.$$
By Lemma \ref{huber}, the assumptions that $\log p=O(n)$ and $H\asymp n^{1/(1+\varepsilon/2)}\to\infty$, we have$$P\Bigg(II>\frac{H}{6}\Bigg)=o(1).$$ By the Cauchy-Schwarz inequality, $|III_t|\leq \sqrt{II}\sqrt{I_t}$. We then get that
$$
\aligned
P\Bigg(|III_t|>\frac{H}{6}\Bigg)\leq &P\Bigg(I_t>\frac{H}{6}\Bigg)+P\Bigg(II>\frac{H}{6}\Bigg)\\
=&o(1).
\endaligned$$
Combining the results above yields that
\begin{equation}\label{part1_egamma}
P\Bigg(\Big|\frac{\|\mathbf{y}_t-\widehat{\boldsymbol{\mu}}\|_2^2}{p}-E(\xi_t^2)\Big|>\frac{H}{2}\Bigg)=o(1).
\end{equation}
Using \eqref{huber_xy_ineq}, we have
\begin{equation}\label{huber_ineq}\aligned
&\frac{1}{n}\sum_{t=1}^{n}\Psi_H\Bigg(\frac{\|\mathbf{y}_t-\widehat{\boldsymbol{\mu}}\|_2^2}{p}-\frac{H}{2}\Bigg)\\
\leq& \frac{1}{n}\sum_{t=1}^{n}\Psi_H\Bigg(\frac{\|\mathbf{y}_t-\widehat{\boldsymbol{\mu}}\|_2^2}{p}-E(\xi_t^2)\Bigg)\\
&-\Big(\frac{H}{2}-E(\xi_t^2)\Big)\frac{1}{n}\sum_{t=1}^{n}\mathbf{1}_{\Big\{\frac{\|\mathbf{y}_t-\widehat{\boldsymbol{\mu}}\|_2^2}{p}-E(\xi_t^2)<H \,\text{and}\, \frac{\|\mathbf{y}_t-\widehat{\boldsymbol{\mu}}\|_2^2}{p}>-\frac{H}{2}\Big\}}.
\endaligned
\end{equation}
By \eqref{part1_egamma}, we have
$$
P\Bigg(\frac{\|\mathbf{y}_t-\widehat{\boldsymbol{\mu}}\|_2^2}{p}-E(\xi_t^2)<H\Bigg)\geq 1-P\Bigg(\Bigg|\frac{\|\mathbf{y}_t-\widehat{\boldsymbol{\mu}}\|_2^2}{p}-E(\xi_t^2)\Bigg|\geq \frac{H}{2}\Bigg)\to 1,
$$
and
$$
P\Bigg(\frac{\|\mathbf{y}_t-\widehat{\boldsymbol{\mu}}\|_2^2}{p}-E(\xi_t^2)>-\frac{H}{2}\Bigg)\geq 1-P\Bigg(\Bigg|\frac{\|\mathbf{y}_t-\widehat{\boldsymbol{\mu}}\|_2^2}{p}-E(\xi_t^2)\Bigg|\geq \frac{H}{2}\Bigg)\to 1.
$$
Therefore, \begin{equation}\label{eva_h2}
P\Bigg(
\Big(\frac{H}{2}-E(\xi_t^2)\Big)\frac{1}{n}\sum_{t=1}^{n}\mathbf{1}_{\Big\{\frac{\|\mathbf{y}_t-\widehat{\boldsymbol{\mu}}\|_2^2}{p}-E(\xi_t^2)<H \,\text{and}\,\frac{\|\mathbf{y}_t-\widehat{\boldsymbol{\mu}}\|_2^2}{p}>-\frac{H}{2}\Big\}}>\frac{H}{4}\Bigg)\to 1. \end{equation}
Combining  \eqref{Psi_Hbound}, \eqref{huber_ineq} and \eqref{eva_h2}  yields that
$$
P\Bigg(\frac{1}{n}\sum_{t=1}^{n}\Psi_H\Big(\frac{\|\mathbf{y}_t-\widehat{\boldsymbol{\mu}}\|_2^2}{p}-\frac{H}{2}\Big)<-\frac{H}{5}\Bigg)\to 1.
$$
Similarly, one can show that
$$
P\Bigg(\frac{1}{n}\sum_{t=1}^{n}\Psi_H\Big(\frac{\|\mathbf{y}_t-\widehat{\boldsymbol{\mu}}\|_2^2}{p}\Big)>\frac{H}{5}\Bigg)\to 1.
$$
Using the definition of $\widehat{E(\xi_t^2)}$ and the continuity of the function $\Psi_H(\cdot)$, we then get that
\begin{equation}\label{hatvbound}
P\Bigg(0< \widehat{E(\xi_t^2)}<\frac{H}{2}\Bigg)\to 1.
\end{equation}
Combining \eqref{Egammat}, \eqref{part1_egamma} and \eqref{hatvbound} yields that
$$E(\gamma_t)=1-o(1). $$}
By Markov's inequality and that $\tilde{p}_H\leq 1$, we get the desired result \eqref{Prob_Hbound}.
$\qed$

\subsection{Proof of Theorems \ref{thm2} and \ref{thm3}}

The bounds in \eqref{gpoet_sigma0} follow from Theorem \ref{prop:pilot} and Theorem 2.1 of \cite{fan2018large}.
The bounds in \eqref{gpoet_sigma1} follow from \eqref{gpoet_sigma0} and Proposition \ref{Consist:v}.$\qed$


%
%


\section*{Supplementary Materials}
 Supplement to ``{\bf Sub-Gaussian High-dimensional covariance matrix estimation under elliptical factor model with $2+\varepsilon$th  moment.''}

{This supplement contains the proofs of Propositions \ref{lem_denominator_xi}-\ref{Consist:v}, and  Theorems \ref{prop:pilot}--\ref{thm3}. }




\bibliographystyle{apalike}
\bibliography{Cov_est_elliptical-aos_template}       



\newpage
\setcounter{page}{1}
\clearpage

\begin{center}
    {\Large \bf Supplement to ``Sub-Gaussian High-Dimensional Covariance Matrix Estimation under Elliptical Factor Model with $2+\varepsilon$th  Moment''}
    
    \bigskip
  
  {\large Yi Ding and Xinghua Zheng }

\end{center}

\pagenumbering{arabic}

\appendixtitleon
\appendixtitletocon
\begin{appendices}

\numberwithin{equation}{section}

\setcounter{equation}{0}
\renewcommand{\theequation}{A.\arabic{equation}}
\section{Proofs}\label{Proof}
The theorem/proposition/lemma/equation numbers below refer to the main article (\cite{DZ24}). In the following proofs, $c$ and $C$ denote generic constants which do not depend on $p$ and $n$ and can vary from place to place.

\section*{Proof of Lemma \ref{mu_est}}
Under Assumptions \ref{asump1} and \ref{sigma_bound}, we have $\max_{1\leq i\leq p}(\boldsymbol{\Sigma})_{ii}=O(1)$.
The desired bound follows from Lemma \ref{huber} with $\delta=O(1/\log p)$ and $\varepsilon=1$, Bonferroni's inequality, and
the assumption that $\log(p)=O(n)$. $\qed$

\section*{Proof of Lemma \ref{normal_quadraticratios}} Write the  eigen decomposition of $\boldsymbol{\Sigma}$ as  $\boldsymbol{\Sigma}=\mathbf{\Gamma}\boldsymbol{\Lambda}\mathbf{\Gamma}^T$ with $\boldsymbol{\Lambda}=\diag(\lambda_1, ..., \lambda_p)$. Denote $\mathbf{z}_{\Gamma;t}=\mathbf{\Gamma}^T\mathbf{z}_t=:(z_{\Gamma;1t}, ..., z_{\Gamma;pt})^T$. Because $\mathbf{z}_t\sim N(0, \mathbf{I})$, and $\mathbf{\Gamma}\mathbf{\Gamma}^T=\mathbf{\Gamma}^T\mathbf{\Gamma}=\mathbf{I}$, we have that $\mathbf{z}_\Gamma \sim N(0, \mathbf{I})$, and
\begin{equation}\label{trans_quadratic}
\frac{\mathbf{z}_t^T\boldsymbol{\Sigma}\mathbf{z}_t}{\mathbf{z}^T_t\mathbf{z}_t}=\sum_{i=1}^N \lambda_i\frac{ z_{\Gamma;it}^2}{\|\mathbf{z}_{\Gamma;t}\|_2^2}.
\end{equation}
By symmetry, $$E\Bigg(\frac{\mathbf{z}_t^T\boldsymbol{\Sigma}\mathbf{z}_t}{\mathbf{z}^T_t\mathbf{z}_t}\Bigg)=\frac{1}{p} \tr(\boldsymbol{\Sigma}).$$ Equation \eqref{E_ratio_trace} follows.

Next, by \eqref{trans_quadratic}, we have
$$
\aligned
E\Bigg(\frac{(\mathbf{z}_t^T\boldsymbol{\Sigma}\mathbf{z}_t)^2}{(\mathbf{z}^T_t\mathbf{z}_t)^2}\Bigg)
=&\sum_{i,j=1}^p \lambda_i\lambda_j E\Bigg(\frac{z_{\Gamma;it}^2z_{\Gamma;jt}^2}{\|\mathbf{z}_{\Gamma;t}\|_2^4}\Bigg)\\
\leq& \big(\tr(\boldsymbol{\Sigma}^2)\big)E\Bigg(\frac{z_{\Gamma;1t}^4}{\|\mathbf{z}_{\Gamma;t}\|_2^4}\Bigg)+\big(\tr(\boldsymbol{\Sigma})\big)^2E\Bigg(\frac{z_{\Gamma;1t}^2z_{\Gamma;2 t}^2}{\|\mathbf{z}_{\Gamma;t}\|_2^4}\Bigg)\\
\leq & \frac{\tr(\boldsymbol{\Sigma}^2)}{p}+\frac{\big(\tr(\boldsymbol{\Sigma})\big)^2}{(p-1)p},
\endaligned
$$
where the last inequality holds because
$$
\sum_{1\leq i,j\leq p}z_{\Gamma;it}^2z_{\Gamma;jt}^2=\|\mathbf{z}_{\Gamma;t}\|_2^4,
$$
which implies that
$$p(p-1)E\Bigg(\frac{z_{\Gamma;1t}^2z_{\Gamma;2 t}^2}{\|\mathbf{z}_{\Gamma;t}\|_2^4}\Bigg)=\sum_{1\leq i\neq j\leq p} E\Bigg(\frac{z_{\Gamma;it}^2z_{\Gamma;j t}^2}{\|\mathbf{z}_{\Gamma;t}\|_2^4}\Bigg)\leq 1,$$ and
$$pE\Bigg(\frac{z_{\Gamma;1t}^4}{\|\mathbf{z}_{\Gamma;t}\|_2^4}\Bigg)=\sum_{i=1}^p E\Bigg(\frac{z_{\Gamma;it}^4}{\|\mathbf{z}_{\Gamma;t}\|_2^4}\Bigg)\leq 1.$$
Inequality \eqref{E_ratio2_trace} follows.
$\qed$

\section*{Proof of Lemma \ref{gaussian_pilot}} The desired bounds are immediate results of Lemmas A.2, A.3 and A.5 of \cite{DLZ19}.$\qed$

For any vector $\mathbf{x}$, let $(\mathbf{x})_i$ be its $i$th entry of $\mathbf{x}$, and for   any matrix $\mathbf{A}$, let $(\mathbf{A})_{ij}$ be its $(i,j)$th entry.

\section*{Proof of Lemma \ref{lem_denominator}}
We start with \eqref{prob_bound_1}. Using \eqref{Y:model},
we have
\begin{equation}\label{l2_decomp}
\aligned
&\frac{\|\widehat{\mathbf{P}}_I(\mathbf{y}_t-\widehat{\boldsymbol{\mu}})\|_2^2}{p\xi_t^2}\\
= &\frac{\|\widehat{\mathbf{P}}_I({\boldsymbol{\mu}}-\widehat{\boldsymbol{\mu}})\|_2^2}{p\xi_t^2}+\frac{\|\widehat{\mathbf{P}}_I\boldsymbol{\Sigma}_0^{1/2}\mathbf{z}_t\|_2^2}{\|\mathbf{z}_t\|_2^2}+2\frac{({\boldsymbol{\mu}}-\widehat{\boldsymbol{\mu}})^T\widehat{\mathbf{P}}^T_I\widehat{\mathbf{P}}_I\boldsymbol{\Sigma}_0^{1/2}\mathbf{z}_t}{\sqrt{p}\xi_t\|\mathbf{z}_t\|_2}\\
=:&I_t+II_t+III_t.
\endaligned
\end{equation}
By Lemma \ref{mu_est} and the assumption that $\xi_t>c$,
$$
\max_{1\leq t\leq n}|I_t|=O_p\Bigg(\frac{\log p}{n}\Bigg).
$$
Define $\mathbf{\Gamma}_{I}=(\phi_{K+1},..., \phi_{p})$.
About term $II_t$, we have
\begin{equation}\label{hatu_decomp}
\aligned
II_t= &\frac{\|\mathbf{\Gamma}_{I}^T\boldsymbol{\Sigma}_0^{1/2}\mathbf{z}_t\|_2^2}{\|\mathbf{z}_t\|_2^2}+\frac{\|(\widehat{\mathbf{P}}_I-\mathbf{\Gamma}_{I}^T)\boldsymbol{\Sigma}_0^{1/2}\mathbf{z}_t\|_2^2}{\|\mathbf{z}_t\|_2^2}+2\frac{\mathbf{z}_t^T\boldsymbol{\Sigma}_0^{1/2}\mathbf{\Gamma}_{I}(\widehat{\mathbf{P}}_I-\mathbf{\Gamma}_{I}^T)\boldsymbol{\Sigma}_0^{1/2}\mathbf{z}_t}{\|\mathbf{z}_t\|^2_2}.
\endaligned
\end{equation}
We have
$E(\|\mathbf{\Gamma}_{I}^T\boldsymbol{\Sigma}_0^{1/2}\mathbf{z}_t\|_2^2/p)=\sigma_a$, and by Lemma \ref{ZL2012},
\begin{equation}\label{eq_u2}
E\Bigg(\Big(\frac{\|\mathbf{\Gamma}_{I}^T\boldsymbol{\Sigma}_0^{1/2}\mathbf{z}_t\|_2^2}{p}-\sigma_a\Big)^{2s}\Bigg)\leq \frac{c_s}{p^{s}},\text{ for all }s\geq 1.
\end{equation}
{Under the assumption that $\log n=O(\log p)$, there exists $\delta>0$ such that $n\leq p^{\delta}$ for all $n$ large enough. We then set $s$ sufficient large such that $s>\delta+1$ and get }
\begin{equation}\label{hatu_decomp_b1}
\aligned
&P\Bigg(\max_{1\leq t\leq n}\Big|\frac{\|\mathbf{\Gamma}_{I}^T\boldsymbol{\Sigma}_0^{1/2}\mathbf{z}_t\|_2^2}{p}-\sigma_a\Big|>\frac{\sigma_a}{5}\Bigg)\\
\leq& cn E\Bigg(\Big(\frac{\|\mathbf{\Gamma}_{I}^T\boldsymbol{\Sigma}_0^{1/2}\mathbf{z}_t\|_2^2}{p}-\sigma_a\Big)^{2s}\Bigg)\\
=& O\Bigg(\frac{n}{p^{s}}\Bigg)\to 0.
\endaligned
\end{equation}
By Lemma \ref{ZL2012} again,
\begin{equation}\label{hatu_decomp_b1_z}
\aligned
&P\Bigg(\max_{1\leq t\leq n}\Big|\frac{\|\mathbf{z}_t\|_2^2}{p}-1\Big|>\frac{1}{5}\Bigg)\\
\leq& cn E\Bigg(\Big(\frac{\|\mathbf{z}_t\|_2^2}{p}-1\Big)^{2s}\Bigg)\\
=& O\Bigg(\frac{n}{p^{s}}\Bigg)\to 0. \endaligned
\end{equation}
By Assumption \ref{sigma_bound} and Lemma \ref{Chatterjee2013}, we have $E\Big(\max_{1\leq t\leq n, 1\leq i\leq p}|(\boldsymbol{\Sigma}_0^{1/2}\mathbf{z}_t)_{i}|\Big)
=O(\sqrt{\log (pn)})$. Hence by the assumptions that $\log n=O(\log p)$ and  $(\log (p))^{2+\gamma}=o(n)$ for some $\gamma>0$, we have
$$\aligned
&P\Big(\max_{1\leq t\leq n, 1\leq i\leq p}(\boldsymbol{\Sigma}_0^{1/2}\mathbf{z}_t)_{i}\geq (\log (p))^{(1+\gamma)/2}\Big)\\
\leq& \frac{E\Big(\max_{1\leq t\leq n, 1\leq i\leq p}|(\boldsymbol{\Sigma}_0^{1/2}\mathbf{z}_t)_{i}|\Big)}{(\log (p))^{(1+\gamma)/2}}\\
=&O\Bigg(\frac{(\log p+\log {n})^{1/2}}{(\log (p))^{(1+\gamma)/2}}\Bigg)=O\Bigg(\frac{1}{(\log (p))^{\gamma/2}}\Bigg)\to 0.
\endaligned
$$
In addition, by Theorem 4.1 of \cite{fan2018large}\footnote{Note that the Spatial Kendall's tau is independent of the scalar process $(\xi_t)$, see P.1402 in \cite{fan2018large}.}, we get that
$\|\widehat{\boldsymbol{\Gamma}_{K}}_{ED}-\boldsymbol{\Gamma}_K\|_2=O_p(\sqrt{(\log p)/n})$. It follows that
\begin{equation}\label{hatu_error}
\|\mathbf{\Gamma}_{I}^T-\widehat{\mathbf{P}}_I\|_2=O_p\Bigg(\sqrt{\frac{\log p}{n}}\Bigg).
\end{equation}
Therefore,
\begin{equation}\label{hatu_decomp_b2}
\aligned
&\max_{1\leq t\leq n}\frac{\|(\widehat{\mathbf{P}}_I-\mathbf{\Gamma}_{I}^T)\boldsymbol{\Sigma}_0^{1/2}\mathbf{z}_t\|_2^2}{p}
\\
\leq& \|(\widehat{\mathbf{P}}_I-\mathbf{\Gamma}_{I}^T)\|_2^2\max_{1\leq t\leq n,1\leq i\leq p}|(\boldsymbol{\Sigma}_0^{1/2}\mathbf{z}_t)_{i}^2|\\
=&O_p\Bigg(\frac{(\log p)^{2+\gamma}}{n}\Bigg)=o_p(1),
\endaligned
\end{equation}
where the last equality holds by the assumption that $(\log (p))^{2+\gamma}=o(n)$.
By \eqref{hatu_decomp_b1} and \eqref{hatu_decomp_b2}, we get
\begin{equation}\label{hatu_decomp_b3}
\aligned
&\max_{1\leq t\leq n}\frac{1}{p}\Big|\mathbf{z}_t^T\boldsymbol{\Sigma}_0^{1/2}\mathbf{\Gamma}_{I}(\widehat{\mathbf{P}}_I-\mathbf{\Gamma}_{I}^T)\boldsymbol{\Sigma}_0^{1/2}\mathbf{z}_t\Big|\\
\leq &\frac{1}{p}\Bigg(\max_{1\leq t\leq n}\|(\widehat{\mathbf{P}}_I-\mathbf{\Gamma}_{I}^T)\boldsymbol{\Sigma}_0^{1/2}\mathbf{z}_t\|_2\Bigg)\Bigg(\max_{1\leq t\leq n}\|\mathbf{\Gamma}_{I}^T\boldsymbol{\Sigma}_0^{1/2}\mathbf{z}_t\|_2\Bigg)\\
=&O_p\Bigg(\frac{(\log p)^{1+\gamma/2}}{\sqrt{n}}\Bigg)=o_p(1).
\endaligned
\end{equation}
Combining \eqref{hatu_decomp}, \eqref{hatu_decomp_b1}, \eqref{hatu_decomp_b1_z}, \eqref{hatu_decomp_b2} and \eqref{hatu_decomp_b3} yields
$$P\Bigg(\max_{1\leq t\leq n}|II_t-\sigma_a|>\frac{\sigma_a}{3}\Bigg)=o(1).$$
About term $III_t$, we  have
$$\max_{1\leq t\leq n}|III_t|\leq \sqrt{\Big(\max_{1\leq t\leq n}I_t\Big)\Big(\max_{1\leq t\leq n}II_t\Big)}=o_p(1). $$
The desired result \eqref{prob_bound_1} follows.

Next, we show \eqref{sum_bound_1}. We have
\begin{equation}\label{flip}
\aligned
&
\frac{1}{n}\sum_{t=1}^n\Bigg(\frac{p\xi_t^2}{\|\widehat{\mathbf{P}}_I(\mathbf{y}_t-\widehat{\boldsymbol{\mu}})\|_2^2}-\frac{1}{\sigma_a}\Bigg)^2\\
\leq& \frac{1}{\sigma_a^2 \min_{1\leq t\leq n}\frac{\|\widehat{\mathbf{P}}_I(\mathbf{y}_t-\widehat{\boldsymbol{\mu}})\|_2^4}{p^2\xi_t^4}}\sum_{t=1}^n\frac{1}{n}\Big(\frac{\|\widehat{\mathbf{P}}_I(\mathbf{y}_t-\widehat{\boldsymbol{\mu}})\|_2^2}{p\xi_t^2}-\sigma_a\Big)^2\\
\leq& \frac{2}{\sigma_a^2 \min_{1\leq t\leq n}\frac{\|\widehat{\mathbf{P}}_I(\mathbf{y}_t-\widehat{\boldsymbol{\mu}})\|_2^4}{p^2\xi_t^4}}\Bigg(\frac{1}{n}\sum_{t=1}^n\Big(\frac{\|\mathbf{\Gamma}_{I}^T\boldsymbol{\Sigma}_0^{1/2}\mathbf{z}_t\|_2^2}{\|\mathbf{z}_t\|_2^2}-\sigma_a\Big)^2\\
&\quad\quad\quad\quad\quad\quad\quad\quad\quad\quad\quad\quad+\frac{1}{n}\sum_{t=1}^n\Big(\frac{\|\mathbf{\Gamma}_{I}^T\boldsymbol{\Sigma}_0^{1/2}\mathbf{z}_t\|_2^2}{\|\mathbf{z}_t\|_2^2}-\frac{\|\widehat{\mathbf{P}}_I(\mathbf{y}_t-\widehat{\boldsymbol{\mu}})\|_2^2}{p\xi_t^2}\Big)^2\Bigg)\\
=:&\frac{2}{\sigma_a^2\min_{1\leq t\leq n}\frac{\|\widehat{\mathbf{P}}_I(\mathbf{y}_t-\widehat{\boldsymbol{\mu}})\|_2^4}{p^2\xi_t^4}}(I+II).
\endaligned
\end{equation}
About term $I$, we have
$$\aligned
I\leq& \frac{2}{n}\sum_{t=1}^n\Big(\frac{\|\mathbf{\Gamma}_{I}^T\boldsymbol{\Sigma}_0^{1/2}\mathbf{z}_t\|_2^2}{p}-\sigma_a\Big)^2+\frac{2}{n}\sum_{t=1}^n\frac{\|\mathbf{\Gamma}_{I}^T\boldsymbol{\Sigma}_0^{1/2}\mathbf{z}_t\|_2^2}{\|\mathbf{z}_t\|_2^2} \Big(\frac{\|\mathbf{z}_t\|_2^2}{p}-1\Big)^2\\
\leq& \frac{2}{n}\sum_{t=1}^n\Big(\frac{\|\mathbf{\Gamma}_{I}^T\boldsymbol{\Sigma}_0^{1/2}\mathbf{z}_t\|_2^2}{p}-\sigma_a\Big)^2+\frac{2C}{n}\sum_{t=1}^n\Big(\frac{\|\mathbf{z}_t\|_2^2}{p}-1\Big)^2,
\endaligned
$$
where the last inequality holds because $$\frac{\|\mathbf{\Gamma}_{I}^T\boldsymbol{\Sigma}_0^{1/2}\mathbf{z}_t\|_2^2}{\|\mathbf{z}_t\|_2^2}\leq \|\mathbf{\Gamma}_{I}^T\boldsymbol{\Sigma}_0\mathbf{\Gamma}_{I}\|_2=O(1).$$
By Jensen's inequality and \eqref{eq_u2}, for all $k\geq 1$,
$$
\aligned
&E\Bigg(\Bigg(\frac{1}{n}\sum_{t=1}^n\Big(
\frac{\|\mathbf{\Gamma}_{I}^T\boldsymbol{\Sigma}_0^{1/2}\mathbf{z}_t\|_2^2}{p}-\sigma_a\Big)^2\Bigg)^{k}\Bigg)\\
\leq& \frac{1}{n}\sum_{t=1}^nE\Bigg(\Big(\frac{\|\mathbf{\Gamma}_{I}^T\boldsymbol{\Sigma}_0^{1/2}\mathbf{z}_t\|_2^2}{p}-\sigma_a\Big)^{2k}\Bigg)\\
=&O\Bigg(\frac{1}{p^k}\Bigg).
\endaligned$$
Hence,
$$
\aligned
P\Bigg(\frac{1}{n}\sum_{t=1}^n\Big(\frac{\|\mathbf{\Gamma}_{I}^T\boldsymbol{\Sigma}_0^{1/2}\mathbf{z}_t\|_2^2}{p}-\sigma_a\Big)^2\geq \frac{\log p}{p}\Bigg)\leq \frac{c}{\log p}\to 0.
\endaligned
$$
Similarly,
\begin{equation}\label{z2_p}
\aligned
P\Bigg(\frac{1}{n}\sum_{t=1}^n\Big(\frac{\|\mathbf{z}_t\|_2^2}{p}-1\Big)^2\geq \frac{\log p}{p}\Bigg)\leq \frac{c}{\log p}\to 0.
\endaligned
\end{equation}
It follows that
\begin{equation}\label{term1}
I=O_p\Bigg(\frac{\log p}{p}\Bigg).
\end{equation}
About term $II$, by \eqref{l2_decomp} and \eqref{hatu_decomp}, and the inequality that $(\sum_{k=1}^K x_k)^2\leq 2^{K-1}\sum_{k=1}^Kx_k^2$, we have
\begin{equation}\label{term2_part2}
\aligned
&\Bigg(\frac{\|\mathbf{\Gamma}_{I}^T\boldsymbol{\Sigma}_0^{1/2}\mathbf{z}_t\|_2^2}{\|\mathbf{z}_t\|_2^2}-\frac{\|\widehat{\mathbf{P}}_I(\mathbf{y}_t-\widehat{\boldsymbol{\mu}})\|_2^2}{p\xi_t^2}\Bigg)^2\\
\leq &
8\frac{\|\widehat{\mathbf{P}}_I({\boldsymbol{\mu}}-\widehat{\boldsymbol{\mu}})\|_2^4}{p^2\xi_t^4}\\
&+8\frac{\|\widehat{\mathbf{P}}_I-\mathbf{\Gamma}_{I}^T\|_2^4(\mathbf{z}_t^T\boldsymbol{\Sigma}_0\mathbf{z}_t)^2}{\|\mathbf{z}_t\|_2^4}\\
&+32\frac{\|{\boldsymbol{\mu}}-\widehat{\boldsymbol{\mu}}\|_2^4\cdot\|\widehat{\mathbf{P}}_I-\mathbf{\Gamma}_{I}^T\|_2^4(\mathbf{z}_t^T\boldsymbol{\Sigma}_0\mathbf{z}_t)^2}{p^2\xi^4_t\|\mathbf{z}_t\|_2^4}\\
&+32\frac{\|{\boldsymbol{\mu}}-\widehat{\boldsymbol{\mu}}\|_2^4(\mathbf{z}_t^T\mathbf{\Gamma}_{I}^T\boldsymbol{\Sigma}_0\mathbf{\Gamma}_{I}\mathbf{z}_t)^2}{p^2\xi^4_t\|\mathbf{z}_t\|_2^4}.
\endaligned
\end{equation}
By Lemma \ref{ZL2012}, we have that
$$E\Bigg(\frac{\mathbf{z}_t^T\boldsymbol{\Sigma}_0\mathbf{z}_t}{p}\Bigg)=1,\quad\text{ and }\quad E\Bigg(\Big(\frac{\mathbf{z}_t^T\boldsymbol{\Sigma}_0\mathbf{z}_t}{p}\Big)^2\Bigg)=O(1).$$ Hence,
\begin{equation}\label{term2_part3}
\frac{1}{n}\sum_{t=1}^n \frac{\mathbf{z}_t^T\boldsymbol{\Sigma}_0\mathbf{z}_t}{p}=O_p(1),\quad\text{ and }\quad\frac{1}{n}\sum_{t=1}^n \Bigg(\frac{\mathbf{z}_t^T\boldsymbol{\Sigma}_0\mathbf{z}_t}{p}\Bigg)^2=O_p(1).
\end{equation}
Because $E(\mathbf{z}_t^T\mathbf{\Gamma}_{I}^T\boldsymbol{\Sigma}_0\mathbf{\Gamma}_{I}\mathbf{z}_t/p)=\sigma_a=O(1)$, we also have
\begin{equation}\label{term2_part4}
\frac{1}{n}\sum_{t=1}^n\frac{ \mathbf{z}_t^T\mathbf{\Gamma}_{I}^T\boldsymbol{\Sigma}_0\mathbf{\Gamma}_{I}\mathbf{z}_t}{p}=O_p(1).
\end{equation}
Combining Lemma \ref{mu_est}, \eqref{hatu_decomp_b1_z}, \eqref{hatu_error}, \eqref{term2_part2}, \eqref{term2_part3} and \eqref{term2_part4} yields that
\begin{equation}\label{term2}
II=O_p\Bigg(\frac{\log p}{n}\Bigg).
\end{equation}
The desired bound \eqref{sum_bound_1} follows from \eqref{prob_bound_1}, \eqref{flip},  \eqref{term1} and \eqref{term2}.

Finally, we show \eqref{trace_x}. Recall that $$\widehat{\mathbf{X}}_t=\frac{\sqrt{p}}{\|\widehat{\mathbf{P}}_I(\mathbf{y}_t-\widehat{\boldsymbol{\mu}})\|_2}(\mathbf{y}_t-\widehat{\boldsymbol{\mu}})=\frac{\sqrt{p}}{\|\widehat{\mathbf{P}}_I(\mathbf{y}_t-\widehat{\boldsymbol{\mu}})\|_2}\Bigg(\boldsymbol{\mu}-\widehat{\boldsymbol{\mu}}+\frac{\sqrt{p}{\xi_t}}{\|\mathbf{z}_t\|_2}\widetilde{\mathbf{X}}_t\Bigg).$$ We have
$$\aligned
\widehat{\mathbf{X}}_t\widehat{\mathbf{X}}_t^T=&\frac{p}{\|\widehat{\mathbf{P}}_I(\mathbf{y}_t-\widehat{\boldsymbol{\mu}})\|_2^2}(\boldsymbol{\mu}-\widehat{\boldsymbol{\mu}})(\boldsymbol{\mu}-\widehat{\boldsymbol{\mu}})^T\\
&+\frac{p^{3/2}\xi_t}{\|\widehat{\mathbf{P}}_I(\mathbf{y}_t-\widehat{\boldsymbol{\mu}})\|_2^2\|\mathbf{z}_t\|_2}(\boldsymbol{\mu}-\widehat{\boldsymbol{\mu}})\widetilde{\mathbf{X}}_t^T\\
&+\frac{p^{3/2}\xi_t}{\|\widehat{\mathbf{P}}_I(\mathbf{y}_t-\widehat{\boldsymbol{\mu}})\|_2^2\|\mathbf{z}_t\|_2}\widetilde{\mathbf{X}}_t(\boldsymbol{\mu}-\widehat{\boldsymbol{\mu}})^T\\
&+\Big(\frac{p^2\xi^2_t}{\|\widehat{\mathbf{P}}_I(\mathbf{y}_t-\widehat{\boldsymbol{\mu}})\|_2^2\|\mathbf{z}_t\|_2^2}-\frac{1}{\sigma_a}\Big)\widetilde{\mathbf{X}}_t\widetilde{\mathbf{X}}_t^T\\
&+\frac{1}{\sigma_a}\widetilde{\mathbf{X}}_t\widetilde{\mathbf{X}}_t^T.
\endaligned$$
Therefore, for any $1\leq i,j\leq p$,
$$\aligned
\Big(\frac{1}{n}\sum_{t=1}^n\widehat{\mathbf{X}}_t\widehat{\mathbf{X}}_t^T\Big)_{ij}=\frac{1}{\sigma_a}(\widetilde{\mathbf{S}}_0)_{ij}&+\frac{1}{n}\sum_{t=1}^n\frac{p}{\|\widehat{\mathbf{P}}_I(\mathbf{y}_t-\widehat{\boldsymbol{\mu}})\|_2^2}(\boldsymbol{\mu}-\widehat{\boldsymbol{\mu}})_i(\boldsymbol{\mu}-\widehat{\boldsymbol{\mu}})_j\\
&+\frac{1}{n}\sum_{t=1}^n\frac{p^{3/2}\xi_t}{\|\widehat{\mathbf{P}}_I(\mathbf{y}_t-\widehat{\boldsymbol{\mu}})\|_2^2\|\mathbf{z}_t\|_2}(\boldsymbol{\mu}-\widehat{\boldsymbol{\mu}})_i(\widetilde{\mathbf{X}}_t)_j\\
&+\frac{1}{n}\sum_{t=1}^n\frac{p^{3/2}\xi_t}{\|\widehat{\mathbf{P}}_I(\mathbf{y}_t-\widehat{\boldsymbol{\mu}})\|_2^2\|\mathbf{z}_t\|_2}(\widetilde{\mathbf{X}}_t)_i(\boldsymbol{\mu}-\widehat{\boldsymbol{\mu}})_j\\
&+\frac{1}{n}\sum_{t=1}^n\Big(\frac{p^2\xi^2_t}{\|\widehat{\mathbf{P}}_I(\mathbf{y}_t-\widehat{\boldsymbol{\mu}})\|_2^2\|\mathbf{z}_t\|_2^2}-\frac{1}{\sigma_a}\Big)(\widetilde{\mathbf{X}}_t)_i(\widetilde{\mathbf{X}}_t)_j\\
=:\frac{1}{\sigma_a}(\widetilde{\mathbf{S}}_0)_{ij}&+I_{ij}+II_{ij}+III_{ij}+IV_{ij}.
\endaligned$$

About term $I_{ij}$, by Lemma \ref{mu_est}, \eqref{prob_bound_1}  and the assumption that $\xi_t>c$, we have
$$
\max_{1\leq i,j\leq p}|I_{ij}|=O_p\Bigg(\frac{\log p}{n}\Bigg).
$$
About terms $II_{ij}$ and $III_{ij}$, we have
$$
\aligned
&\max_{1\leq i,j\leq p}\big(|II_{ij}|,|III_{ij}|\big)\\
\leq &\Bigg(\max_{1\leq t\leq n}\frac{p^{3/2}\xi_t}{\|\widehat{\mathbf{P}}_I(\mathbf{y}_t-\widehat{\boldsymbol{\mu}})\|_2^2\|\mathbf{z}_t\|_2}\Bigg)\cdot\Bigg(\max_{1\leq i,j\leq p}\frac{1}{n}\sum_{t=1}^n|(\widetilde{\mathbf{X}}_t)_i|\cdot|(\boldsymbol{\mu}-\widehat{\boldsymbol{\mu}})_j|\Bigg)\\
\leq &
\Bigg(\max_{1\leq t\leq n}\frac{p^{3/2}\xi_t}{\|\widehat{\mathbf{P}}_I(\mathbf{y}_t-\widehat{\boldsymbol{\mu}})\|_2^2\|\mathbf{z}_t\|_2}\Bigg)\cdot\sqrt{\max_{1\leq i\leq p}\frac{1}{n}\sum_{t=1}^n\big((\widetilde{\mathbf{X}}_t)_i\big)^2}\cdot\|\boldsymbol{\mu}-\widehat{\boldsymbol{\mu}}\|_{\max}.
\endaligned
$$
It follows from  \eqref{Gaussian:Sig}, \eqref{prob_bound_1} and
\eqref{hatu_decomp_b1_z} that
$$
\max_{1\leq i,j\leq p}\big(|II_{ij}|,|III_{ij}|\big)=O_p\Bigg(\sqrt{\frac{\log p}{n}}\Bigg).
$$
About term $IV_{ij}$, we have
$$\aligned
|IV_{ij}|\leq &\frac{1}{n}\sum_{t=1}^n\frac{p}{\|\mathbf{z}_t\|_2^2}\Big|\frac{p\xi^2_t}{\|\widehat{\mathbf{P}}_I(\mathbf{y}_t-\widehat{\boldsymbol{\mu}})\|_2^2}-\frac{1}{\sigma_a}\Big|\cdot |\widetilde{x}_{it}\widetilde{x}_{jt}|\\
&+\frac{1}{n}\sum_{t=1}^n\frac{1}{\sigma_a} \Bigg|\frac{p}{\|\mathbf{z}_t\|_2^2}-1\Bigg||\widetilde{x}_{it}\widetilde{x}_{jt}|,
\endaligned
$$
where $\widetilde{x}_{it}=(\widetilde{\mathbf{X}}_t)_i$.
By the Cauchy-Schwarz inequality,
we have
$$\aligned
\max_{1\leq i,j\leq p}|IV_{ij}|\leq &\Bigg(\max_{1\leq t\leq n}\frac{p}{\|\mathbf{z}_t\|_2^2}\sqrt{\frac{1}{n}\sum_{t=1}^n\Big|\frac{p\xi^2_t}{\|\widehat{\mathbf{P}}_I(\mathbf{y}_t-\widehat{\boldsymbol{\mu}})\|_2^2}-\frac{1}{\sigma_a}\Big|^2}\Bigg)\cdot \sqrt{\max_{1\leq i\leq p}\frac{1}{n}\sum_{t=1}^n\widetilde{x}_{it}^4}\\
&+\frac{1}{\sigma_a}\sqrt{\frac{1}{n}\sum_{t=1}^n\Big(\frac{p}{\|\mathbf{z}_t\|_2^2}-1\Big)^2}\sqrt{\max_{1\leq i\leq p}\frac{1}{n}\sum_{t=1}^n \widetilde{x}^4_{it}}.
\endaligned
$$
By Bernstein's inequality, we have $\max_{1\leq i\leq p}(\sum_{t=1}^n\widetilde{x}_{it}^4/n)=O_p(1)$. By \eqref{sum_bound_1}, \eqref{hatu_decomp_b1_z} and \eqref{z2_p}, we then get
$$
\max_{1\leq i\leq p}|IV_{ij}|=O_p\Bigg(\sqrt{\frac{\log p}{n}}+\sqrt{\frac{\log p}{p}}\Bigg).
$$
Combining the results above yields
\begin{equation}\label{max_error_hatx}
\Big\|\Big(\frac{1}{n}\sum_{t=1}^n\widehat{\mathbf{X}}_t\widehat{\mathbf{X}}_t^T\Big)-\frac{1}{\sigma_a}\widetilde{\mathbf{S}}_0\Big\|_{\max}=O_p\Bigg(\sqrt{\frac{\log p}{n}}+\sqrt{\frac{\log p}{p}}\Bigg).
\end{equation}
In addition, by \eqref{Gaussian:Sig}, $$\frac{\tr(\widetilde{\mathbf{S}}_0)}{p}=\frac{\tr(\boldsymbol{\Sigma}_0)}{p} +O_{p}\Bigg(\sqrt{\frac{\log p}{n}}\Bigg)=1+O_{p}\Bigg(\sqrt{\frac{\log p}{n}}\Bigg).$$
The desired result \eqref{trace_x} follows.
$\qed$

\section*{ Proof of Lemma \ref{prop:pilot:Sigma}}
The desired bound \eqref{pilot:Sigma_0:Sig} follows from \eqref{trace_x} and \eqref{max_error_hatx}.
The bound \eqref{pilot:Sigma_0:Lam0} then follows from \eqref{pilot:Sigma_0:Sig}, Assumption~\ref{assump_factor_structure}, Lemma \ref{gaussian_pilot} and Weyl's Theorem.

Finally, about \eqref{pilot:Sigma_0:Gam0}, by Assumption \ref{assump_factor_structure}, \eqref{pilot:Sigma_0:Sig}, \eqref{pilot:Sigma_0:Lam0} and the sin$\theta$ Theorem,
\begin{equation}\label{gam_l2}
\|\widehat{\boldsymbol{\Gamma}}_{K}-\widetilde{\boldsymbol{\Gamma}}_{K}\|_2=O_p\Bigg(\sqrt{\frac{\log p}{n}}+\sqrt{\frac{\log p}{p}}\Bigg).
\end{equation}
 By the definitions of $\widehat{\boldsymbol{\Gamma}}_{K}$ and $\widetilde{\boldsymbol{\Gamma}}_{K}$, we have $\widehat{\phi}_{k}=(1/\widehat{\lambda}_{0;k})\widehat{\boldsymbol{\Sigma}}_{0}\widehat{\phi}_{k}$ and $\widetilde{\phi}_{k}=(1/\widetilde{\lambda}_{0;k})\widetilde{\mathbf{S}}_{0}\widetilde{\phi}_{k}$. Hence, for all $1\leq k\leq K$,
\begin{equation}\label{pilot_decompo_gamma}
\aligned
\|\widehat{\phi}_{k}-\widetilde{\phi}_{k}\|_{\max}\leq &\Bigg|1-\frac{\widehat{\lambda}_{0;k}}{\widetilde{{\lambda}}_{0;k}}\Bigg|\cdot\|\widehat{\phi}_{k}\|_{\max}\\
&+\frac{1}{\widetilde{\lambda}_{0;k}}\|(\widehat{\boldsymbol{\Sigma}}_{0}-\widetilde{\mathbf{S}}_{0})\widehat{\phi}_{k}\|_{\max}\\
&+\frac{1}{\widetilde{\lambda}_{0;k}}\|\widetilde{\mathbf{S}}_{0}(\widehat{\phi}_{k}-\widetilde{\phi}_{k})\|_{\max}.
\endaligned
\end{equation}
Note that $$\|(\widehat{\boldsymbol{\Sigma}}_{0}-\widetilde{\mathbf{S}}_{0})\widehat{\phi}_{k}\|_{\max}\leq\sqrt{p} \|\widehat{\boldsymbol{\Sigma}}_{0}-\widetilde{\mathbf{S}}_{0}\|_{\max}\|\widehat{\phi}_{k}\|_2,$$
and
$$\|\widetilde{\mathbf{S}}_{0}(\widehat{\phi}_{k}-\widetilde{\phi}_{k})\|_{\max}\leq \sqrt{p}\|\widetilde{\mathbf{S}}_{0}\|_{\max}\|\widehat{\phi}_{k}-\widetilde{\phi}_{k}\|_2.$$
By Assumption \ref{assump_factor_structure}, Lemma \ref{gaussian_pilot}, \eqref{pilot:Sigma_0:Sig} and \eqref{pilot:Sigma_0:Lam0}, we have $\|\widehat{\boldsymbol{\Sigma}}_0\|_{\max}=O_p(1)$, $\|\widehat{\boldsymbol{\Sigma}}_0\|_{\max}=O_p(1)$, and there exists $c>0$ such that with probability approaching one, for all large $p$,
\begin{equation}\label{lowerbound_lambda}
\widehat{\lambda}_{0;K}>cp, \quad\text{ and} \quad \widetilde{\lambda}_{0;K}>cp.
\end{equation}
Because $\widehat{\boldsymbol{\Gamma}}_{K}=\widehat{\boldsymbol{\Sigma}}_0\widehat{\boldsymbol{\Gamma}}_{K}\widehat{\boldsymbol{\Lambda}}_{0K}^{-1}$, it follows that
\begin{equation}\label{gamhat_max}
\|\widehat{\boldsymbol{\Gamma}}_{K}\|_{\max}\leq \sqrt{p} \|\widehat{\boldsymbol{\Sigma}}_0\|_{\max}\|\widehat{\boldsymbol{\Gamma}}_{K}\|_2\|\widehat{\boldsymbol{\Lambda}}_{0K}^{-1}\|_{2}=O_p\Bigg(\frac{1}{\sqrt{p}}\Bigg).
\end{equation}
The desired bound \eqref{pilot:Sigma_0:Gam0} follows from \eqref{pilot:Sigma_0:Sig}, \eqref{pilot:Sigma_0:Lam0} and \eqref{gam_l2}--\eqref{gamhat_max}.
$\qed$

\end{appendices}

\end{document}